\documentclass[10pt, reqno, a4paper]{amsart}
\usepackage{amsmath,amsthm,amssymb,amsfonts}

\usepackage{enumitem}
\usepackage{a4wide}
\usepackage{bm}
\usepackage{mathtools}
\usepackage{thmtools}
\usepackage[dvipsnames]{xcolor}
\usepackage[colorlinks=true]{hyperref}
\usepackage{graphicx}
\usepackage{orcidlink}

\theoremstyle{plain}
\declaretheorem[name=Proposition, numberwithin=section]{proposition}
\declaretheorem[name=Theorem, sibling=proposition]{theorem}
\declaretheorem[name=Lemma, sibling=proposition]{lemma}
\declaretheorem[name=Corollary, sibling=proposition]{corollary}

\theoremstyle{definition}
\declaretheorem[name=Remark, sibling=proposition]{remark}

\DeclareMathOperator{\E}{\mathsf{E}}
\renewcommand{\P}{\mathsf{P}}

\DeclareMathOperator{\N}{\mathsf{N}}
\DeclareMathOperator{\R}{\mathsf{R}}
\DeclareMathOperator{\D}{\mathsf{D}}
\DeclareMathOperator{\F}{\mathcal{F}}
\renewcommand{\d}{\mathrm{d}}
\DeclareMathOperator{\eqd}{\,\overset{\mathsf{D}}{=}\,}

\newcommand{\numberthis}{\addtocounter{equation}{1}\tag{\theequation}}

\begin{document}
\numberwithin{equation}{section}

\title[Inference for SDEs driven by Hermite processes]{Inference for SDEs driven by Hermite processes}

\author{Petr \v Coupek\orcidlink{0000-0002-5360-6095}}
\address{
	Charles University, 
	Faculty of Mathematics and Physics, 
	Sokolovsk\' a 83, 
	Prague 8, 186 75, 
	Czech Republic. }
\email{coupek@karlin.mff.cuni.cz}

\author{Pavel K\v r\'i\v z\orcidlink{0000-0001-5773-9541}}
\address{
	Charles University, 
	Faculty of Mathematics and Physics, 
	Sokolovsk\' a 83, 
	Prague 8, 186 75, 
	Czech Republic.}
\email{kriz@karlin.mff.cuni.cz}

\begin{abstract}
In the paper, we address parametric and non-parametric estimation for nonlinear stochastic differential equations with additive Hermite noise with possibly nonlinear scaling. We assume that a single trajectory of the solution is observed discretely and we propose estimators of the Hurst parameter and the Hermite order of the driving process as well as of the average noise intensity and noise intensity function. The estimators are based on the weighted quadratic variation whose properties are used, in particular, to prove weak consistency of the proposed estimators under in-fill asymptotics.
\end{abstract}

\keywords{Hermite process, Estimation, Hurst parameter, Hermite order, Noise intensity}

\subjclass[2020]{60G22; 62M09}



\maketitle


\section{Introduction}

In the article, we focus on estimation for Hermite processes (HPs) and stochastic differential equations (SDEs) driven by them. The family of HPs arises naturally as the limit of suitably normalized sums of strongly correlated random variables in what is sometimes called the non-central limit theorem; see, e.g., \cite{DobMaj79}. These processes exhibit self-similarity, stationarity of increments as well as long-range dependence, and they can be described as multiple Wiener--It\^o integrals with suitable deterministic kernels; see, e.g., \cite{Tud13}. The family of HPs is parametrized by two parameters---the Hurst parameter $H\in (1/2,1)$, that governs the self-similarity of the process, and the Hermite order $q\in\N$, that, roughly speaking, measures the degree of non-Gaussianity of the process (more precisely, this is the order of the Wiener chaos in which the process lives, see section \ref{sec:HP}). The most well-known examples of HPs is the fractional Brownian motion with Hurst parameter $H>1/2$, which is the HP of order $1$, and the Rosenblatt process, which is the HP of order $2$. We refer to, e.g., \cite{DecUst99, ManNes68} for the definition and basic properties of fractional Brownian motions and to \cite{BiaHuOksZha08} and the references therein for a more thorough exposition. We also refer to \cite{Taq11} for an accessible exposition of Rosenblatt processes, to \cite{Tud08} for their thorough analysis, and to, e.g., \cite{AbrPip06, Alb98, CouDunDun22, DawLoo22, DavKer23, GarTorTud12, KerNouSakVii21, Pip04} for some of their finer properties. Of course, the general class of HPs has also been considered in the literature and we refer to, e.g., \cite{PipTaq17, Tud13, Tud23} for an overview and to, e.g.,  \cite{Aya20, AyaHamLoo24, BaiTaq14, Loo23, MorOod86} for some finer results and extensions.

As far as statistical inference is concerned, much of the literature has been devoted to inference for models driven by fractional Brownian motions; see, e.g., the monographs \cite{KubMisRal2017, PraRao10} and the references therein. On the other hand, the literature on statistical inference for models with general Hermite noise is still quite limited. 

Of course, one basic task is the estimation of the Hurst parameter $H$. An estimator of $H$ for a Rosenblatt process in a high-frequency regime based on quadratic variations is given in \cite{TudVie09} and this approach is generalized to Hermite processes in \cite{ChrTudVie09, ChrTudVie11}. The estimator was shown not to be asymptotically normal. Strongly consistent asymptotically normal estimators were recently proposed in \cite{AyaTud24} (based on modified quadratic variations) and in \cite{LooTud25} (based on wavelet variations). 

On the other hand, apart from testing Gaussianity, the problem of identification of parameter $q$ seems not to have been addressed in the literature so far. This is rather surprising both in view of the multitude of situations in which Gaussianity of the driving noise is not a reasonable assumption; see, e.g., \cite{ChaRoi20,CheEtAl19, Dom15,SkaMabSch13,SkaMabSch14,WitGotLanFraGreHei19}, and in the fact that some of the proposed estimators of $H$ depend on the properties of the driving HP (e.g.\ the estimator in \cite{ChrTudVie11} depends on $q$).

As far as inference for SDEs driven by Hermite noise is concerned, apart from estimation of parameters $H$ and $q$, there are other parameters that need to be inferred; namely, the drift and noise intensity coefficients. In \cite{AssTud20}, an estimator of the Hurst parameter in a linear SDE with additive Hermite noise (the Hermite Ornstein--Uhlenbeck process) based on quadratic variation is proposed and analyzed and in \cite{NouTra19}, a moment-type estimator of the drift coefficient is studied. Estimation of all the parameters (drift, noise intensity, and Hurst parameter) for SDEs with additive Rosenblatt noise has been recently treated in \cite{CouKriMas25}. 

In fact, the present article can be viewed as a continuation and significant extension of the results from \cite{CouKriMas25}. Namely, we consider the nonlinear SDE driven by additive Hermite noise with nonlinear scaling
	\[ 
		\d{Y}_t = g(Y_t)\d{t} + f(t)\d{Z}_t^{H,q}
	\] 
on a fixed time horizon $[0,T]$ where $g:\R\to\R$ is a deterministic nonlinear function that is locally Lipschitz and satisfies a certain Lyapunov-type condition, $f:[0,T]\to \R$ is a deterministic bounded function, and $Z^{H,q}$ is an HP with Hurst parameter $H\in (1/2,1)$ and Hermite order $q\in\N$. 

The above equation is a generalization of the problem considered in \cite{CouKriMas25} in the sense that here we allow for general HP and we also allow for a non-constant noise intensity (in \cite{CouKriMas25}, only Rosenblatt noise with constant noise intensity was considered). 

We assume that we observe a single trajectory of the solution $Y$ at discrete points and we aim to obtain estimators of the Hurst parameter $H$, the noise intensity function $f$, and, to the best of our knowledge for the very first time, of the Hermite order $q$. 

As in \cite{ChrTudVie11}, we base our estimators on the normalized quadratic variation. It turns out that in the non-Gaussian case, the behavior of the estimators computed from the solution reduces to the behavior of the estimators of the noise term. This is not surprising as the solution $Y$ can be viewed as a (regular) permutation of the (irregular) noise term. Our estimator of the Hurst parameter is shown to be weakly consistent under the in-fill asymptotics and, moreover, its asymptotic distribution is non-Gaussian (in fact, regardless of the Hermite order of the driving HP, the limit belongs to the second Wiener chaos). It is robust and it does not depend on any of the other parameters. A similar result on the asymptotics is also obtained for the estimator of the $L^2$-average noise intensity. Additionally, we propose a non-parametric estimator of the noise intensity function and prove that it is weakly consistent with respect to the $L^2$-distance. These two estimators are plug-in in nature---for their computation, one has to estimate the Hurst parameter from the data first. Finally, we propose a procedure to identify the Hermite order of the driving process. We consider the situation where the noise intensity is constant and known (e.g.\ if the driving process is the standard HP of an unknown order) and our estimator is shown to be weakly consistent under the in-fill asymptotics. If the noise intensity is unknown, we can consistently estimate the Hermite order if there is no drift in the equation.

\subsection*{Organization of the article} The article is organized as follows. 
In section \ref{sec:prelim}, we define the family of HPs and the Wiener integral with respect to them. In section \ref{sec:qvar_of_wiener_integral}, we investigate the properties of the quadratic variation of the Wiener integral. In particular, we prove that the integral has a (suitably normalized) quadratic variation (\autoref{thm: vars joint convergence}) and we identify its speed of convergence and asymptotic distribution (\autoref{thm: 2-vars distrib limit}). Section \ref{sec:estimation_wiener_integral} is devoted to the estimation for the Wiener integral. Initially, we propose an estimator of the Hurst parameter of the driving HP and identify its (non-Gaussian) asymptotic distribution (\autoref{cor: H estimator}). Subsequently, we give an estimator of the $L^2$-average of the noise intensity function,  identify its asymptotic distribution (\autoref{cor: avg intensity estim}), and then we propose a nonparametric estimator of the noise intensity function and prove its weak consistency with respect to the $L^2$-distance (\autoref{thm: weak consistency of noise estimator}). Finally, we propose an estimator of the Hermite order for a scaled HP with known scaling and prove its weak consistency (\autoref{thm: est H' and q}). Section \ref{sec:SDEs} is devoted to the extension of the results to SDEs and the article is concluded with a discussion on the estimation of the Hermite order of the scaled HP in the situation when the scaling factor is unknown in \autoref{app:unknown_scaling} and with some auxiliary lemmas in \autoref{app:aux_lemmas}.

\subsection*{Conventions} Throughout the article, the following notation is used. 
\begin{itemize}
	\itemsep0em
	\item For a function $A: [0,T]\to\R$, we write $A_{s,t} := A_t-A_s$ for $s,t\in [0,T]$, $s\leq t$.
	\item We write $A\lesssim B$ if there exists a finite positive constant $C$ such that $A\leq C B$.
	\item For a set $A$, we use the symbol $\# A$ to denote the number of its elements.
	\item We write $A\propto B$ if there exists a finite positive constant $C$ such that $A=CB$. 
	\item We write $a_n\asymp b_n$ as $n\to\infty$ for two sequences $\{a_n\}_n, \{b_n\}_n\subset (0,\infty)$ if there exists a finite positive constant $C$ such that the inequality $b_n/C \leq a_b \leq Cb_n$ holds for all $n$.
	\item We use the symbols $\|\cdot\|_{L^p}$ for the $L^p$-norm, $\|\cdot\|_\infty$ for the supremum norm, and $\|\cdot\|_{C^\alpha}$ for the $\alpha$-H\"older norm.
\end{itemize}

\section{Preliminaries}
\label{sec:prelim}

Let $T\in (0,\infty)$ and let $(\Omega,\F,\P)$ be a complete probability space rich enough to carry a standard Wiener process $W=(W_t)_{t\in [0,T]}$. We assume that $\F$ is generated by $W$. For $q\in\N_0$, let $I_q$ denote the Wiener--It\^o multiple integral of order $q$ with respect to $W$ (cf., e.g., \cite[Section 1.1.2]{Nua06}) and let $\mathcal{H}_q$ denote the corresponding $q$\textsuperscript{th} Wiener chaos (cf., e.g., \cite[Section 1.1.1]{Nua06}).

\subsection{Definition of Hermite processes}
\label{sec:HP}
For the analysis contained in this article, we consider the following finite-time interval representation of HP; see, e.g., \cite[Definition 7]{NouNuaTud10}, \cite[Theorem 1.1]{PipTaq10}, or \cite[Theorem 3.1]{Tud13}. For $H\in (1/2,1)$, $q\in\N$, and $t\in [0,T]$, define $L^{H,q}_t: [0,T]^q\to\R$ by
	\begin{equation}
	\label{eq:LHq}
		L^{H,q}_t(x) 
			:= c_{H,q} \left[\int_{x_1\vee \ldots\vee x_q}^t \prod_{i = 1}^q \left(\frac{u}{x_i}\right)^{\frac{1}{2}-\frac{1-H}{q}} (u-x_i)^{-\left(\frac{1}{2}+\frac{1-H}{q}\right)}\d{u}\right] \bm{1}_{(0,t)^q}(x)
	\end{equation}
for $x=(x_1,\ldots,x_q)\in [0,T]^q$ where 
	\begin{equation}
	\label{eq:cHq}
	c_{H,q} := \left(\frac{H(2H-1)}{q!\mathrm{B}\left(\frac{1}{2}-\frac{1-H}{q},\frac{2(1-H)}{q}\right)^q}\right)^\frac{1}{2}
	\end{equation}
in which $\mathrm{B}(\cdot\,,\cdot)$ is the Beta function. We have $L_t^{H,q}\in L^2([0,T]^q)$ and we can therefore define
	\[
		Z_t^{H,q} := I_q(L^{H,q}_t), \quad t\in [0,T].
	\]
The process $Z^{H,q}=(Z_t^{H,q})_{t\in [0,T]}$ is called \emph{the Hermite process of order $q$ and Hurst index $H$.} The value of the normalizing constant $c_{H,q}$ is chosen in such a way that process $Z^{H,q}$ is normalized, i.e. that $\E (Z_1^{H,q})^2=1$ holds. It is a standard calculation (see, e.g., \cite{Tud13}) to show that process $Z^{H,q}$ is centered and its covariance function is given by 
	\begin{equation}
	\label{eq:ACF}
	\E Z_s^{H,q}Z_t^{H,q} = \frac{1}{2}\left(s^{2H} + t^{2H} - |t-s|^{2H}\right), \quad s,t\in [0,T].
	\end{equation}

\subsection{Wiener integral with respect to a Hermite process}
In this section, we recall the definition of a Wiener integral with respect to the HP from, e.g., \cite{CouMasOnd22} and some of its properties.

Let $Z^{H,q}$ be the HP of order $q\in\N$ with Hurst index $H\in(1/2,1)$. Denote by $\mathcal{E}(0,T)$, the space of deterministic step functions, i.e. each $f\in\mathcal{E}(0,T)$ is a function $f:[0,T]\to\R$ that satisfies
	\begin{equation}
	\label{eq:elem_f}
		f = \sum_{i=0}^{n-1}f^i\bm{1}_{[t_i, t_{i+1})}
	\end{equation}
for some $\{f^i\}_{i}\subset\R$ and $\{0=t_0\leq t_1\leq\ldots \leq t_{n}=T\}$. For such a step function, let us set 
	\[ 
		i(f) := \sum_{i=0}^{n-1} f^iZ_{t_i,t_{i+1}}^{H,q}.
	\]
Consider also the operator $L^{H,q}: \mathcal{E}(0,T)\to L^2([0,T]^q)$ defined by 
	\[ 
		(L^{H,q}f)(x) := c_{H,q} \int_{x_1\vee\ldots\vee x_q}^T f(u) \prod_{i=1}^q \left(\frac{u}{x_i}\right)^{\frac{1}{2}-\frac{1-H}{q}} (u-x_i)_+^{-\left(\frac{1}{2} + \frac{1-H}{q}\right)}\d{u}
	\]
for $f\in \mathcal{E}(0,T)$ and $x=(x_1,x_2,\ldots, x_q)\in [0,T]^q$ where $c_{H,q}$ is defined by \eqref{eq:cHq}. Note that the map $\langle\cdot,\cdot\rangle_{\mathcal{H}}: \mathcal{E}(0,T)\times\mathcal{E}(0,T)\to\R$ defined by 
	\[ 
		\langle f,g\rangle_{\mathcal{D}_H(0,T)} := H(2H-1)\int_0^T \int_0^T f(u)g(v)|u-v|^{2H-2}\d{u}\d{v}\]
for $f,g\in \mathcal{E}(0,T)$, is an inner product on $\mathcal{E}(0,T)$ and denote by $\|\cdot\|_{\mathcal{D}_H(0,T)}$ the induced norm. For $f\in\mathcal{E}(0,T)$, we immediately obtain that 
	\[ i(f) = I_q(L^{H,q}f) \]
holds and 
	\[ \|i(f) \|_{L^2(\Omega)}^2 = q! \| L^{H,q}f\|_{L^2([0,T]^q)}^2 = \|f\|_{\mathcal{D}_H(0,T)}^2 \]
where the first equality follows by the isometry property of multiple Wiener--It\^o integrals, cf., e.g., \cite[point (iii) on p.\ 9]{Nua06}, while the second equality follows by a straightforward calculation using Fubini's theorem and the equality 
	\begin{equation}
	\label{eq:beta}
		(uv)^{\alpha} \int_0^{u\wedge v} x^{-2\alpha} (u-x)^{\alpha-1}(v-x)^{\alpha-1}\d{x} = \mathrm{B}(\alpha,1-2\alpha) |u-v|^{2\alpha-1}, \quad u\neq v,
	\end{equation}
that holds for $\alpha\in (0,1/2)$. In this manner, we obtain that the operator $i$ is a linear isometry between the space $\mathcal{E}(0,T)$ endowed with the norm $\|\cdot\|_{\mathcal{D}_H(0,T)}$ and the linear span of $\{i(f): f\in\mathcal{E}(0,T)\}$ endowed with the norm $\|\cdot\|_{L^2(\Omega)}$; see, e.g., \cite{MaeTud07} for a similar computation. As such, it can be uniquely extended to a linear isometry between the completion $\mathcal{D}_H(0,T)$ of $\mathcal{E}(0,T)$ with respect to $\|\cdot\|_{\mathcal{D}_H(0,T)}$ and the closure of the the linear span of $\{i(f): f\in\mathcal{E}(0,T)\}$ in the norm $L^2(\Omega)$. This extension is denoted by $\int_0^T (\cdot)\d{Z}_s^{H,q}$ and it can be shown that its natural domain $\mathcal{D}^H(0,T)$ is the homogeneous Sobolev--Slobodeckij space $\dot{W}^{\frac{1}{2}-H,2}(0,T)$, cf. \cite[Proposition 2.6]{CouMasOnd22}. Because of the embedding $L^\frac{1}{H}(0,T)\subseteq \dot{W}^{\frac{1}{2}-H,2}(0,T)$, we immediately obtain the estimate in the following lemma that we will need throughout the article:

\begin{proposition}
\label{lem:bddness_of_int}
For every $H\in (1/2,1)$, $q\in\N$, and $f\in L^\frac{1}{H}(0,T)$, we have
	\[
		\left\| \int_0^T f_s\d{Z}_s^{H,q}\right\|_{L^2(\Omega)} \lesssim \| f\|_{L^\frac{1}{H}(0,T)}.
	\]
\end{proposition}

\section{Quadratic variation of the integral with respect to Hermite processes}
\label{sec:qvar_of_wiener_integral}
Fix $T\in (0,\infty)$ and consider the sequence of equidistant partitions of the interval $[0,T]$ and its sub-sequence that are defined for $n\in\N$ by
	\begin{equation}
	\label{eq:dyadic_partitions}
		\left\{T_j^n:=T\frac{j}{K_n}\right\}_{j=0}^{K_n} \subseteq \left \{t_i^n = T \frac{i}{n}\right\}_{i=0}^{n} 
	\end{equation}
where $K_n\in \N$. We assume that that $n/K_n\to \infty$ as $n\to\infty$. For each $n\in\N$, denote the mesh sizes of $\{t_i^n\}_i$ and $\{T_j^n\}_j$ by $\delta_n$ and $\Delta_n$, respectively, i.e. 
	\[ 
		\delta_n := t_{i+1}^n-t_{i}^n= \frac{T}{n}, \qquad \Delta_n := T_{j+1}^n - T_{j}^n = \frac{T}{K_n},
	\]	
and denote also the number of points from the original partition $\{t_i^n\}_i$ contained in one sub-partition interval $[T_j^n,T_{j+1}^n]$ by $D_n$, i.e.
	\[D_n := \# \{i: t_i^n \in [T_j^n,T_{j+1}^n) \}.
	\]
Denote also, for a process $X=(X_t)_{t\in [0,T]}$, two points $s,t\in [0,T]$ satisfying $s\leq t$, and $n\in\N$,
\[
V_{n,[s,t]}^{2}(X) := \sum_{i: t_i^n \in [s,t)}|X_{t_i^n,t_{i+1}^n}|^{2}.  
\]

\subsection{Convergence of quadratic variation} 
In what follows, we show that the quadratic variation of the Wiener integral with respect to the HP converges when suitably normalized and we also identify its limit. The main result of this section is the following

\begin{theorem}
\label{thm: vars joint convergence}
Let $f\in L^2(0,T)$ be a given deterministic function and let $X=(X_t)_{t\in [0,T]}$ be defined by $X_t:= \int_0^t f_s\d{Z}_s^{H,q}$, $t\in [0,T]$. Then for every $p\in [1,\infty)$, there is the convergence
\begin{equation}
	\label{eq: Unif p-var}
 \sum_{j=0}^{K_n-1}\left| \frac{1}{\delta_n^{2H-1}}V_{n,[T_j^n,T_{j+1}^n]}^{2}(X) - \int_{T_j^n}^{T_{j+1}^n} \left|f_s\right|^{2}\d{s}\right| \quad\xrightarrow[n \to \infty]{L^p(\Omega)}\quad0.
\end{equation}
\end{theorem}

\begin{proof}
Consider a sequence of functions $\{f^m\}_{m\in\N}$ such that 
	\begin{equation*}
        \|f-f^m\|_{L^{2}(0,T)} \xrightarrow[m \to \infty]{} 0
    \end{equation*}
and such that for each $m\in\mathbb{N}$, function $f^m$ satisfies
     \begin{equation*}
        f^m = \sum_{j=0}^{\kappa_m-1} f^{m,j}\bm{1}_{[s_j^m,s_{j+1}^m)},
    \end{equation*}
for some $\kappa_m\in\mathbb{N}$, $\{f^{m,j}\}_{j=0}^{\kappa_m-1}\subset \R$, and a partition $\{s_j^m\}_{j=0}^{\kappa_m}$ of the interval $[0,T]$ such that for each fixed $n\in\N$ we either have $\{ T_j^n\}_{j=0}^{K_n}\subseteq \{s_j^m\}_{j=0}^{\kappa_m}$ or $\{ T_j^n\}_{j=0}^{K_n}\supseteq \{s_j^m\}_{j=0}^{\kappa_m}$. Consider also a sequence of processes $\{X^m\}_{m\in\mathbb{N}}$ where for each $m\in\N$, process $X^m=(X_t^m)_{t\in [0,T]}$ is defined by 
	\[
		X^m_t := \int_0^t f_s^m \d Z_s^{H,q}, \quad t\in [0,T].
	\]
	
For $m,n\in\mathbb{N}$ fixed, apply the triangle inequality to obtain the estimate
    \begin{align*}
      \sum_{j=0}^{K_n-1}\left| \frac{1}{\delta_n^{2H-1}} V_{n,[T_j^n,T_{j+1}^n]}^{2}(X) - \int_{T_j^n}^{T_{j+1}^n} \left|f_s\right|^{2}\d{s} \right|  &  \\
        & \hspace{-7cm} \leq  \sum_{j=0}^{K_n-1} \frac{1}{\delta_n^{2H-1}} \left| V_{n,[T_j^n,T_{j+1}^n]}^{2}(X) - V_{n,[T_j^n,T_{j+1}^n]}^{2}(X^m) \right| \\
        & \hspace{-6cm} + \sum_{j=0}^{K_n-1}\left| \frac{1}{\delta_n^{2H-1}} V_{n,[T_j^n,T_{j+1}^n]}^{2}(X^m) - \int_{T_j^n}^{T_{j+1}^n} \left|f^m_s\right|^{2}\d{s} \right|   \\
        & \hspace{-6cm} + \sum_{j=0}^{K_n-1} \left|\int_{T_j^n}^{T_{j+1}^n} (|f_s^m|^2-|f_s|^2)\d{s}\right| \\
        & \hspace{-7cm} =: a_n^m + b_n^m + c^m.\\
    \end{align*}
We aim to show that all the three terms in the expression above converge to zero in the $L^p(\Omega)$-norm. We will prove this in multiple steps.
    
    \textit{Step 1.} \quad In the first step, we aim to show that 
    	\begin{equation}
    	\label{eq:step_1}
    		\sup_{n\in\mathbb{N}}\|a_n^m\|_{L^p(\Omega)}\xrightarrow[m\to\infty]{} 0.
    	\end{equation}
   For $m, n\in\mathbb{N}$, we have the estimate  
	\begin{align*}
      \left\| a_n^m  \right\|_{L^{p}(\Omega)}  
      	& = \left\|\sum_{j=0}^{K_n-1} \frac{1}{\delta_n^{2H-1}} \left| V_{n,[T_j^n,T_{j+1}^n]}^{2}(X) - V_{n,[T_j^n,T_{j+1}^n]}^{2}(X^m) \right|\right\|_{L^p(\Omega)} \\
		& = \left\|\sum_{j=0}^{K_n-1} \frac{1}{\delta_n^{2H-1}}  \left|\sum_{i: t_i^n \in [T_j^n,T_{j+1}^n)}|X_{t_i^n,t_{i+1}^n}|^{2}  - \sum_{i: t_i^n \in [T_j^n,T_{j+1}^n)}|X_{t_i^n,t_{i+1}^n}^m|^{2} \right|\right\|_{L^p(\Omega)}\\
		& \leq  \frac{1}{\delta_n^{2H-1}} \sum_{i=0}^{n-1} \left\| |X_{t_i^n,t_{i+1}^n}|^{2}  - |X_{t_i^n,t_{i+1}^n}^m |^{2} \right\|_{L^p(\Omega)}
	\end{align*}
As both random variables $X_{t_i^n,t_{i+1}^n}$ and $X^m_{t_i^n,t_{i+1}^n}$ belong to the $q$\textsuperscript{th} Wiener chaos $\mathcal{H}_q$, their squares belong to $\bigoplus_{k=0}^{2q}\mathcal{H}_k$ by the product formula for multiple Wiener--It\^o integrals; see, e.g., \cite[Proposition 1.1.2]{Nua06}. As the space $\bigoplus_{k=0}^{2q}\mathcal{H}_k$ is linear, the difference of the squares also belongs to $\bigoplus_{k=0}^{2q}\mathcal{H}_k$ and therefore, it has equivalent moments, i.e.\ its $L^p(\Omega)$-norm can be estimated by its $L^{1}(\Omega)$-norm; see, e.g., \cite[Proposition 2.2]{CouMasOnd22}. Thus we obtain the estimate
	\[ 
	\left\| a_n^m  \right\|_{L^{p}(\Omega)}   \lesssim \frac{1}{\delta_n^{2H-1}} \sum_{i=0}^{n-1} \left\| |X_{t_i^n,t_{i+1}^n}|^{2}  - |X_{t_i^n,t_{i+1}^n}^m|^{2} \right\|_{L^1(\Omega)}.
	\]
Now we can proceed similarly as in the proof of \cite[Lemma 4.2]{GueNua05} with $1/H$ replaced by $2$ and obtain
		\begin{equation*}
		\left\| a_n^m  \right\|_{L^{p}(\Omega)} \lesssim \frac{1}{\delta_n^{2H-1}} \left(\E V_{n,[0,T]}^2(X-X^m)\right)^\frac{1}{2} \left( \left(\E V_{n,[0,T]}^2(X)\right)^{\frac{1}{2}} +\left(\E V_{n,[0,T]}^2(X^m)\right)^{\frac{1}{2}}\right).
	  \end{equation*}
The subsequent application of \autoref{lem: pVar bounded by Lp} yields
		\begin{equation*}
		\left\| a_n^m  \right\|_{L^{p}(\Omega)} \lesssim \left\| f - f^m  \right\|_{L^2(0,T)} ( \left\| f \right\|_{L^2(0,T)} +\left\| f^m \right\|_{L^2(0,T)})
	  \end{equation*}
   from which convergence \eqref{eq:step_1} follows.
		
    \textit{Step 2.}\quad In the second step, we aim to show that for $m\in\mathbb N$, we have    
    \[ 
    	\| b_n^m\|_{L^p(\Omega)} \xrightarrow[n\to\infty]{} 0.
    \]
Let $m\in\mathbb{N}$ and for $n\in\mathbb{N}$ define a new partition $\{\tilde{T}_j^n\}_{j=1}^{\tilde{K}_n}$ as the finer of the two partitions $\{T_j^n\}_{j=0}^{K_n}$ and $\{s_j^m \}_{j=0}^{\kappa_m}$. Notice that we can write
    \begin{align*} 
    	b_n^m =& \sum_{j=0}^{K_n-1}\left| \frac{1}{\delta_n^{2H-1}} V_{n,[T_j^n,T_{j+1}^n]}^{2}(X^m) -  \int_{T_j^n}^{T_{j+1}^n} \left|f^m_s\right|^{2}\d{s} \right| \\
    	 \leq &  \sum_{j=0}^{\tilde{K}_n-1} \left| \frac{1}{\delta_n^{2H-1}} V_{n,[\tilde{T}_j^n,\tilde{T}_{j+1}^n]}^{2}(X^m) -  \int_{\tilde{T}_j^n}^{\tilde{T}_{j+1}^n} \left|f^m_s\right|^{2}\d{s} \right| \\
    = & \sum_{j=0}^{\kappa_m-1} \sum_{k \in S_j(n)} \left| \frac{1}{\delta_n^{2H-1}} V_{n,[\tilde{T}_k^n,\tilde{T}_{k+1}^n]}^{2}(X^m) - |f^{m,j} |^2 (\tilde{T}_{k+1}^n- \tilde{T}_k^n ) \right|,
    \end{align*}
    where $S_j(n) := \{k : [\tilde{T}_k^n,\tilde{T}_{k+1}^n] \subseteq [s_j^m, s_{j+1}^m ]\}$.  Further denote 
	\[\tilde{D}_n := \# \{i: t_i^n \in [\tilde{T}_k^n,\tilde{T}_{k+1}^n) \},
	 \quad  \tilde{\Delta}_n := \tilde{T}_{k+1}^n - \tilde{T}_{k}^n.
	\]
    and recall that neither $\tilde{D}_n$ nor $\tilde{\Delta}_n$ depend on $k$ and that $\tilde{D}_n \xrightarrow[n\to\infty]{} \infty$. The use of stationarity of increments together with $H$-self-similarity of $Z^{H,q}$ yields 
    \begin{align*}
     \left\| b_n^m  \right\|_{L^{p}(\Omega)} & = \left\| \sum_{j=0}^{\kappa_m-1} \sum_{k \in S_j(n)} \left| \frac{1}{\delta_n^{2H-1}} V_{n,[\tilde{T}_k^n,\tilde{T}_{k+1}^n]}^{2}(X^m) - |f^{m,j} |^2 (\tilde{T}_{k+1}^n- \tilde{T}_k^n ) \right| \right\|_{L^{p}(\Omega)}\\
   		& \leq  \sum_{j=0}^{\kappa_m-1} \sum_{k \in S_j(n)} \left\|  \frac{1}{\delta_n^{2H-1}} V_{n,[\tilde{T}_k^n,\tilde{T}_{k+1}^n]}^{2}(X^m) - |f^{m,j} |^2 (\tilde{T}_{k+1}^n - \tilde{T}_k^n)  \right\|_{L^{p}(\Omega)}\\ 
		& = 	
		\sum_{j=0}^{\kappa_m-1} \sum_{k \in S_j(n)} \Bigg\|   \frac{1}{\delta_n^{2H-1}} \sum_{i : t_i^n \in [\tilde{T}_k^n,\tilde{T}_{k+1}^n)} |f^{m,j} Z^{H,q}_{t_i^n,t_{i+1}^n}|^{2}
		 - |f^{m,j}|^2  (\tilde{T}_{k+1}^n - \tilde{T}_k^n) \Bigg\|_{L^{p}(\Omega)} \\
		& = \sum_{j=0}^{\kappa_m-1} |f^{m,j}|^2 \sum_{k \in S_j(n)}   \left\| \frac{1}{\delta_n^{2H-1}} \sum_{i=0}^{\tilde{D}_n-1}  \delta_n^{2H} |Z^{H,q}_{i,i+1}|^{2} -  \tilde{\Delta}_n \right\|_{L^{p}(\Omega)} \\
    	& = \sum_{j=0}^{\kappa_m-1} |f^{m,j}|^2  (s_{j+1}^m - s_j^m)  \left\|  \frac{1}{\tilde{D}_n}\sum_{i=0}^{\tilde{D}_n-1}|Z^{H,q}_{i,i+1}|^{2} - 1 \right\|_{L^{p}(\Omega)}  \xrightarrow[n \to \infty]{} 0, \\
	\end{align*}
where the $L^p(\Omega)$ convergence of the averages on the last line follows by the mixing property of the increments of the HP (see \cite[Theorem 8.3.1]{Sam16}) and the mean (von Neumann's) ergodic theorem.

\textit{Step 3.}\quad  Finally, we aim to show that 
	\begin{equation}
	\label{eq:cm}
		c_m\xrightarrow[m\to\infty]{} 0.
	\end{equation}
To this end, fix $m\in\N$, apply the mean value theorem similarly as in the proof of \cite[Lemma 4.2]{GueNua05}, and use H\"{o}lder's inequality to obtain
\begin{align*}
 c^m & =\sum_{j=0}^{K_n-1}\left|\int_{T_j^n}^{T_{j+1}^n} |f_s^m|^2-|f_s|^2\d{s}\right| \\
 	 & \leq \sum_{j=0}^{K_n-1} \int_{T_j^n}^{T_{j+1}^n} \left| |f_s^m|^2-|f_s|^2\right|\d{s} \\
     & = \int_{0}^{T} \left| |f_s^m|^2-|f_s|^2\right|\d{s} \\
     & \leq  \int_0^T  2 |f_s^m-f_s|\left(|f_s^m|+|f_s|\right) \d{s}\\
     & =  2 \|f^m-f\|_{L^{2}(0,T)}\left(\|f^m\|_{L^{2}(0,T)}+\|f\|_{L^{2}(0,T)}\right)
\end{align*}
from which convergence \eqref{eq:cm} follows. This concludes the proof of \eqref{eq: Unif p-var}.
\end{proof}

\subsection{Speed of convergence and asymptotic distribution of quadratic variation}

In what follows, we analyze the speed of convergence and asymptotic distribution of the quadratic variation of the Wiener integral 
with respect to the HP. To this end, let $f \in L^2(0,T)$ and denote the approximation error of the quadratic variation of the Wiener integral with respect to the HP on the interval $[0,T]$ by 
\[
	\mathcal{V}_n(f) := \frac{1}{\delta_n^{2H-1}} V_{n,[0,T]}^2\left(\int_0^\bullet f_s\d{Z}_s^{H,q}\right) - \int_0^T f_s^2 \d s.
\]

There is the following non-central limit theorem for the asymptotic distribution of the error. The claim is an extension of \cite[Theorem 3.2]{ChrTudVie11} and the proof has a similar structure.

\begin{theorem}
\label{thm: 2-vars distrib limit}
Assume either that $q\geq 2$ holds or, if $q=1$, that we have $H>3/4$. Further let
	\begin{equation}
	\label{eq:H'}
    	H': = 1+ \frac{H-1}{q}
    \end{equation}
and let $f$ be a bounded, piecewise $\alpha$-H\"older continuous function on $[0,T]$ for some $\alpha > 2-2H'$\footnote{That is, there exists a finite partition $\{0=s_0<s_1<\ldots<s_m=T\}$ such that for each $j\in \{0,1,\ldots, m-1\}$, function $f$ restricted to the open subinterval $(s_j,j_{j+1})$ is $\alpha$-H\"older continuous.}. Then there is the convergence
    \[ 
    	\delta_n^{2H'-2}  \mathcal{V}_n(f) \quad  \xrightarrow[n \to \infty]{L^2(\Omega)} \quad d_{H',q}\, \int_0^T f_s^2 \d Z_s^{2H'-1,2},   
    \]
where
    \[
    	d_{H',q} := \frac{2^\frac{1}{2}  q[(2H'-2)q+1][(H'-1)q+1]}{[(4H'-3)(2H'-1)]^\frac{1}{2} [(2H'-2)(q-1)+1][(H'-1)(q-1)+1]}.
    \]
\end{theorem}

\begin{proof}
Define a sequence of step functions $\{f^n\}_n$ by 
    \[
        f^n := \sum_{i=0}^{n-1} f^{n,i} \bm{1}_{[t_i^n,t_{i+1}^n)}, \quad n\in\N,
    \]
with
    \[
        f^{n,i} := \frac{1}{\delta_n}\int_{t_i^n}^{t_{i+1}^n} f_s \d{s}, \quad i\in \{0,1,\ldots, n\}.
    \]
    
We split the argument into several steps. In the first step, we show that the behavior of $\mathcal{V}_n(f)$ for $n\to\infty$ is equivalent to the behavior of $\mathcal{V}_n(f^n)$ for $n\to \infty$. More precisely, we prove the following claim:

	\begin{enumerate}[label = \textit{Claim \arabic*} , leftmargin=*]
		\item \label{claim3} It holds that 
			\[ \delta_n^{2H'-2}(\mathcal{V}_n(f) - \mathcal{V}_n(f^n)) \quad\xrightarrow[n\to\infty]{L^2(\Omega)}\quad 0.\]
	\end{enumerate}

As the second step, we focus on $\mathcal{V}_n(f^n)$ and we find its Wiener chaos expansion. More precisely, we prove the following claim:
	
	\begin{enumerate}[label = \textit{Claim \arabic*} , leftmargin=*]
	\setcounter{enumi}{1}
		\item \label{claim0} Term $\mathcal{V}_n(f^n)$ admits the decomposition 
			\begin{equation}
				\label{eq:decomposition}
				\mathcal{V}_n(f^n) = \mathcal{T}_{2q}^n + c_{2q-2}\mathcal{T}_{2q-2}^n + \ldots + c_4\mathcal{T}_4^n + c_2\mathcal{T}_2^n
			\end{equation}
			where $c_{2q-2k}$ and $\mathcal{T}_{2q-2k}$ are for $k\in \{0,1,\ldots, q-1\}$ defined by 
				\[
				c_{2q-2k} := k! {q\choose k}^2
				\]
			and
				\[
					\mathcal{T}_{2q-2k}^n : = I_{2q-2k} \left( \delta_n^{1-2H} \sum_{i=0}^{n-1} (f^{n,i})^2 L^{H,q}_{t_i^n,t_{i+1}^n}\otimes_k L^{H,q}_{t_i^n,t_{i+1}^n}\right).
				\]
			In the above expressions, symbol $L^{H,q}$ denotes the kernel defined by formula \eqref{eq:LHq} and symbol $\otimes_r$, $r\in \{0,1,\ldots, q\}$, denotes the contraction operator \cite[p.\ 10]{Nua06}.
	\end{enumerate}

As the decomposition \eqref{eq:decomposition} is orthogonal, it follows that in order to assess the behavior of $\mathcal{V}_n(f^n)$ as $n\to\infty$ it suffices to analyze the behavior of the individual terms. In fact, we prove that the behavior is governed by the component in the second Wiener chaos. More precisely, the following claims will conclude the proof of the theorem:

	\begin{enumerate}[label = \textit{Claim \arabic*} , leftmargin=*]
	\setcounter{enumi}{2}
		\item \label{claim1} There is the following convergence
			\[
				\delta_n^{2H' -2} \mathcal{T}_2^n \quad\xrightarrow[n\to\infty]{L^2(\Omega)}\quad \frac{d_{H',q}}{c_2} \int_0^T f_s^2\d{Z}_s^{2H'-1,2}.
			\] 
		\item \label{claim2} For every $k\in \{0,1, \ldots, q-2\}$, there is the convergence
			\[ 
				\delta_n^{2H'-2} \mathcal{T}_{2q-2r}^n \quad\xrightarrow[n\to\infty]{L^2(\Omega)} \quad 0.
			\]
	\end{enumerate}
	
We now prove the individual claims separately.


\textit{Proof of \ref{claim3}.} As the random variable $\mathcal{V}_n(f)-\mathcal{V}_n(f^n)$ belongs to $\bigoplus_{j=0}^{2q}\mathcal{H}_j$ by the product formula for multiple Wiener--It\^o integrals, see, e.g., \cite[Proposition 1.1.3]{Nua06}, it has equivalent moments; see, e.g., \cite[Proposition 2.2]{CouMasOnd22}, and we obtain the estimate
    \begin{equation}
    \label{eq:n-vn-vn_2}
         \delta_n^{2H'-2} \left\| \mathcal{V}_n(f) - \mathcal{V}_n(f^n)\right\|_{L^{2}(\Omega)} \lesssim \delta_n^{2H'-2} \left\| \mathcal{V}_n(f) - \mathcal{V}_n(f^n)\right\|_{L^{1}(\Omega)}.
    \end{equation}
Further, we have 
    \begin{align*}
        \left\| \mathcal{V}_n(f) - \mathcal{V}_n(f^n)\right\|_{L^{1}(\Omega)}  & \\
        & \hspace{-2cm} \leq \delta_n^{1-2H}\E \left| V_{n,[0,T]}^2 \left( \int_0^\bullet f_s \d Z_s^{H,q}\right)^2  -  V_{n,[0,T]}^2 \left( \int_{0}^{\bullet} f^n_s \d Z_s^{H,q}\right)^2 \right| \\
        &\hspace{-1cm} + \left|\int_0^T f_s^2 \d s - \int_0^T (f_s^n)^2 \d s  \right| \\
        & \hspace{-2cm} \lesssim \left\| f-f^n \right\|_{L^2(0,T)} \left(\left\| f \right\|_{L^2(0,T)} + \left\| f^n \right\|_{L^2(0,T)} \right) \numberthis\label{eq:vn-vn_2}
    \end{align*}
by \autoref{lem: difference of pVars}. Moreover, the construction of $f^n$ and H\"older continuity of $f$ guarantee that 
    \[
        \left\| f- f^n \right\|_{L^2(0,T)} \lesssim \delta_n^{\alpha} 
    \]
holds. Inserting this last estimate into \eqref{eq:vn-vn_2} and the thus obtained estimate into \eqref{eq:n-vn-vn_2}, we obtain
    \begin{equation*}
         \delta_n^{2H'-2} \left\| \mathcal{V}_n(f) - \mathcal{V}_n(f^n)\right\|_{L^{2}(\Omega)} \lesssim \delta_n^{\alpha + 2H'-2} \xrightarrow[n \to \infty]{}0. 
    \end{equation*}


\textit{Proof of \ref{claim0}.} Note first that by using the product formula for multiple Wiener--It\^o integrals from \cite[Proposition 1.1.3]{Nua06}, we obtain the expansion
	\[
		\left(Z_{t_i^n,t_{i+1}^n}^{H,q}\right)^2 = I_q\left(L_{t_i^n,t_{i+1}^n}^{H,q}\right)^2 = \sum_{r=0}^q r! {q\choose r}^2 I_{2q-2r}\left( L^{H,q}_{t_i^n,t_{i+1}^n}\otimes_r L^{H,q}_{t_i^n,t_{i+1}^n}\right)
	\]
and it follows that we have	
	\[ 
		\mathcal{V}_n(f^n) = \delta_n^{1-2H}\sum_{i=0}^{n-1} (f^{n,i})^2 \left( \sum_{r=0}^q r! {q\choose r}^2 I_{2q-2r}\left( L^{H,q}_{t_i^n, t_{i+1}^n}\otimes_r L^{H,q}_{t_i^n,t_{i+1}^n}\right) - \delta_n^{2H}\right).
	\]
Direct calculation of the $q$\textsuperscript{th} term of the second sum in the expression above, using equality \eqref{eq:beta}, yields the equality
	\[ 
	q!L^{H,q}_{t_i^n,t_{i+1}^n} \otimes_q L^{H,q}_{t_i^n,t_{i+1}^n} = \delta_n^{2H}
	\] 
and the decomposition in \eqref{eq:decomposition} is obtained.

 
\textit{Proof of \ref{claim1}.} Note first that both $ \delta_n^{2H'-2} \mathcal{T}_2^n$ and $ \frac{d_{H',q}}{c_2} \int_0^T f_s^2\d{Z}_s^{2H'-1,2}$ belong to the second Wiener chaos so that to show the desired convergence, it suffices to prove that kernel of the former  converges to the kernel of the latter in $L^2([0,T]^2)$. In particular, we need to prove the convergence 
	\[ 
		\delta_n^{2H'-2}\delta_n^{1-2H}\sum_{i=0}^{n-1} (f^{n,i})^2 L^{H,q}_{t_i^n,t_{i+1}^n}\otimes_{q-1} L^{H,q}_{t_i^n,t_{i+1}^n}\quad \xrightarrow[n\to\infty]{L^2([0,T]^2)}\quad \frac{d_{H',q}}{c_2} L^{2H'-1,2}f^2.
	\]
We aim to use the dominated convergence theorem. By appealing to equality \eqref{eq:beta}, we obtain
	\begin{multline*}
		 (L^{H,q}_{t_i^n,t_{i+1}^n}\otimes_{q-1} L^{H,q}_{t_i^n,t_{i+1}^n})(x,y) \\
			=  c_{H,q}^2\mathrm{B}(H'-1/2,2-2H')^{q-1} \int_{t_i^n}^{t_{i+1}^n}\int_{t_i^n}^{t_{i+1}^n} k^{H'}(u,x)k^{H'}(v,y)|u-v|^{(2H'-2)(q-1)}\d{u}\d{v}
	\end{multline*}
for almost every $(x,y)\in [0,T]^2$ where
	\[ 
		k^{H'}(u,x) := \left(\frac{u}{x}\right)^{H'-\frac{1}{2}} (u-x)_+^{H'-\frac{3}{2}},
	\]
so that 
	\begin{align*}
		\delta_n^{2H'-2}\delta_n^{1-2H}\sum_{i=0}^{n-1} (f^{n,i})^2 (L^{H,q}_{t_i^n,t_{i+1}^n}\otimes_{q-1} L^{H,q}_{t_i^n,t_{i+1}^n})(x,y) = A_{1}^n(x,y) + A_{2}^n(x,y)
	\end{align*}
holds with 
	\begin{align*}
		A_{1}^n(x,y) & := c_{H,q}^2\mathrm{B}(H'-1/2,2-2H')^{q-1} \delta_n^{2H'-2}\delta_n^{1-2H} \sum_{i=0}^{n-1} (f^{n,i})^2 \\
			&\qquad \int_{t_i^n}^{t_{i+1}^n}\int_{t_i^n}^{t_{i+1}^n} k^{H'}(t_i^n,x)k^{H'}(t_{i}^n ,y)|u-v|^{(2H'-2)(q-1)}\d{u}\d{v}, \\
		A_{2}^n(x,y) & := c_{H,q}^2\mathrm{B}(H'-1/2,2-2H')^{q-1} \delta_n^{2H'-2}\delta_n^{1-2H}\sum_{i=0}^{n-1}(f^{n,i})^2  \\
			&\qquad \int_{t_i^n}^{t_{i+1}^n}\int_{t_i^n}^{t_{i+1}^n} (k^{H'}(u,x)k^{H'}(v,y)-k^{H'}(t_i^n,x)k^{H'}(t_{i}^n ,y))|u-v|^{(2H'-2)(q-1)}\d{u}\d{v}.
	\end{align*}
By straightforward calculation, we obtain 
	\[
		A_{1}^n(x,y) = \frac{2c_{H,q}^2\mathrm{B}(H'-1/2,2-2H')^{q-1}}{((2H'-2)(q-1)+1)((2H'-2)(q-1)+2)} \delta_n \sum_{i=0}^{n-1} (f^{n,i})^2 k^{H'}(t_i^n,x)k^{H'}(t_i^n,y) 
	\]
and, as there is the convergence
   	\[
   		\delta_n \sum_{i=0}^{n-1} (f^{n,i})^2 k^{H'}(t_i^n,x)k^{H'}(t_i^n,y) \xrightarrow[n\to\infty]{} \int_0^T f^2_u k^{H'}(u,x)k^{H'}(u,y)\d{u},
   	\]
we obtain the convergence
	\begin{equation}
	\label{eq:A1n_ptws_conv}
		A_{1}^n(x,y) \xrightarrow[n\to\infty]{} \frac{2c_{H,q}^2 \mathrm{B}(H'-1/2, 2-2H')^{q-1}}{c_{2H'-1,2}((2H'-2)(q-1)+1)((2H'-2)(q-1)+2)} (L^{2H'-1,2}f^2)(x,y).
	\end{equation}
where the constant on the right-hand side is exactly $d_{H',q}/c_2$. Note that for $i,j\in \{0,1, \ldots, n-1\}$, we have that 
	\begin{align*}
		\E Z_{t_i^n,t_{i+1}^n}^{2H'-1}Z_{t_j^n,t_{j+1}^n}^{2H'-1} & \propto \int_{t_i^n}^{t_{i+1}^n}\int_{t_j^n}^{t_{j+1}^n} \left(\int_0^T k^{H'}(u,x) k^{H'}(v,x)\d{x}\right)^2\d{u}\d{v} \\
			& \geq \int_{t_i^n}^{t_{i+1}^n}\int_{t_j^n}^{t_{j+1}^n}\left(\int_0^T k^{H'}(t_i^n,x) k^{H'}(t_j^n,x)\d{x}\right)^2\d{u}\d{v} \\
			& = \delta_n^2 \left(\int_0^T k^{H'}(t_i^n,x) k^{H'}(t_j^n,x)\d{x}\right)^2
	\end{align*}
and that 
	\[ 
		\E Z_{t_i^n,t_{i+1}^n}^{2H'-1}Z_{t_j^n,t_{j+1}^n}^{2H'-1} \propto |t_j^n-t_{i+1}^n|^{4H'-2} + |t_{j+1}^n-t_{i}^n|^{4H'-2} - |t_{i+1}^n-t_{j+1}^n|^{4H'-2} - |t_i^n-t_j^n|^{4H'-2} 
	\] 
from \eqref{eq:ACF} so that 
	\begin{align*}
		\|A_1^n\|_{L^2([0,T]^2)}^2 & \propto \delta_n^2 \sum_{i=0}^{n-1}\sum_{j=0}^{n-1} (f^{n,i})^2(f^{n,j})^2 \int_{t_i^n}^{t_{i+1}^n}\int_{t_j^n}^{t_{j+1}^n} \left(\int_0^T k^{H'}(u,x) k^{H'}(v,x)\d{x}\right)^2\d{u}\d{v} \\
		& \leq \|f\|^4_\infty \delta_n^2 \sum_{i=0}^{n-1}\sum_{j=0}^{n-1} \int_{t_i^n}^{t_{i+1}^n}\int_{t_j^n}^{t_{j+1}^n} \left(\int_0^T k^{H'}(u,x) k^{H'}(v,x)\d{x}\right)^2\d{u}\d{v} \\
		& \lesssim \|f\|^4_\infty\sum_{i=0}^{n-1}\sum_{j=0}^{n-1} \E Z_{t_i^n,t_{i+1}^n}^{2H'-1}Z_{t_j^n,t_{j+1}^n}^{2H'-1} \\
		& \propto \|f\|_\infty^4 \delta_n^{4H'-3} \sum_{i=0}^{n-1}\sum_{j=0}^{n-1} \left(|j-i+1|^{4H'-2} + |j-i-1|^{4H'-2} -2 |i-j|^{4H'-2}\right)\\
		& \propto \|f\|_\infty^4 \delta_n^{4H'-3} \sum_{\ell = -n+1}^{n-1} \left( |\ell+1|^{4H'-2}+|\ell-1|^{4H'-2}-2 |\ell|^{4H'-2}\right).
	\end{align*}
As we have $|\ell+1|^{4H'-2}+|\ell-1|^{4H'-2}-2 |\ell|^{4H'-2}\asymp |\ell|^{4H'-4}$ as $|\ell|\to\infty$, we have that 
	\[
	\sum_{\ell = -n+1}^{n-1} \left( |\ell+1|^{4H'-2}+|\ell-1|^{4H'-2}-2 |\ell|^{4H'-2}\right) \asymp  \sum_{\ell=-n+1}^{n-1} |\ell|^{4H'-4} \asymp \delta_n^{-4H'+3}, \qquad n\to\infty,
	\]
and we obtain that 
	\[
		\|A_1^n\|_{L^2([0,T]^2)}^2 \asymp \delta_n^{4H'-3}\delta_n^{-4H'+3} =1, \qquad n\to\infty,
	\]
or, in other words, that sequence $\{A_1^n\}_n$ is convergent, and therefore bounded, in $L^2([0,T]^2)$. This, together with the pointwise convergence \eqref{eq:A1n_ptws_conv}, yields 
	\[
		A_1^n\quad \xrightarrow[n\to\infty]{L^2([0,T]^2)}\quad \frac{d_{H',q}}{c_2}\, L^{2H'-1,2}f^2.
	\]
Moreover, we have 
  	\begin{multline*}
  		 |A_2^n(x,y)| \lesssim \|f\|_\infty^2 \delta_n^{2H'-2}\delta_n^{1-2H} \\ \sum_{i=0}^{n-1} \int_{t_i^n}^{t_{i+1}^n}\int_{t_i^n}^{t_{i+1}^n} (k^{H'}(u,x)k^{H'}(v,y)-k^{H'}(t_i^n,x)k^{H'}(t_{i}^n ,y))|u-v|^{(2H'-2)(q-1)}\d{u}\d{v}
	\end{multline*}
and as we have that $\|f\|_\infty<\infty$, the convergence \[A_2^n \quad\xrightarrow[n\to\infty]{L^2([0,T]^2)}\quad 0\] now follows as in the proof of \cite[Theorem 3.2]{ChrTudVie11} (see also \cite{TudVie09}).
  
  
\textit{Proof of \ref{claim2}.} By using It\^o's isometry for multiple Wiener--It\^o integrals from, e.g., \cite[point (iii) on p.\ 9]{Nua06}, Fubini's theorem, and formula \eqref{eq:beta}, we obtain in a few steps that there is the estimate 
	\begin{multline*}
		\E (\mathcal{T}_{2q-2k}^n)^2 \lesssim \|f\|_\infty^4 \delta_n^{2-4H} \sum_{i=0}^{n-1}\sum_{j=0}^{n-1} \int_{t_i^n}^{t_{i+1}^n}\int_{t_i^n}^{t_{i+1}^n}\int_{t_j^n}^{t_{j+1}^n}\int_{t_j^n}^{t_{j+1}^n} \\ (|u-u'| |v-v'|)^{(2H'-2)(q-k)}(|u-v||u'-v'|)^{(2H'-2)k}\d{u}\d{v}\d{u'}\d{v'}.
	\end{multline*}
for $k\in \{0,1,\ldots, q-2\}$. As we have that $\|f\|_\infty<\infty$ and the claim follows by the second half of the proof of \cite[Proposition 3.1]{ChrTudVie11}.


\end{proof}

We give two immediate corollaries of the above theorem.

\begin{corollary}
Under the assumptions of \autoref{thm: 2-vars distrib limit}, we have that 
	\[ 
		\delta_n^{4H'-4}\E\mathcal{V}_n(f)^2 \quad\xrightarrow[n\to\infty]{}\quad \tilde{d}_{H',q} \int_0^T\int_0^T f_u^2f_v^2|u-v|^{4H'-4}\d{u}\d{v}
	\]
where 
	\[ 
		\tilde{d}_{H',q} := \frac{2q^2 [(2H'-2)q+1]^2[(H'-1)q+1]^2}{[(2H'-2)(q-1)+1]^2[(2H'-2)(q-1)+1]^2}.
	\]
\end{corollary}

\begin{corollary}
Under the assumptions of \autoref{thm: 2-vars distrib limit}, we have that 
	\[
		\frac{1}{\delta_n^{2H-1}}V_{n,[0,T]}^2\left(\int_0^\bullet f_s\d{Z}_s^{H,q}\right) \quad \xrightarrow[n\to\infty]{a.s.} \quad \int_{0}^{T} f_s^2\d{s}.
	\]
\end{corollary}

\begin{proof}
The claim follows similarly as the in the proof of \cite[Theorem 2.1]{CouKriMas25} by using Markov's inequality, equivalence of moments on a finite Wiener chaos, and a Borel--Cantelli lemma.
\end{proof}

\section{Estimation for Wiener integral with respect to Hermite process}
\label{sec:estimation_wiener_integral}

In this section, we apply the obtained convergence results to parameter estimation for Wiener integrals with respect to Hermite processes. In particular, we consider a trajectory of process $X:=\int_0^\bullet f_s\d{Z}_s^{H,q}$ on the interval $[0,T]$ observed discretely at points $\{t_i^n\}_{i=0}^n$ and we aim to estimate the Hurst index $H$, the noise intensity $|f|$, and we also aim to identify the Hermite order $q$.

\subsection{Estimation of the Hurst index}

Let us begin with an estimator of the Hurst index $H$. In the following result, such an estimator $\widehat{H}_n$ is proposed and its asymptotic distribution is identified.

\begin{corollary}
\label{cor: H estimator}
Assume that $q\geq 2$ or, if $q=1$, that $H>3/4$. Let also $f$ be a non-degenerate ($f\neq 0$), bounded, piecewise $\alpha$-H\"older continuous function on $[0,T]$ for some $\alpha>2-2H'$ where $H'$ is defined by \eqref{eq:H'}. Define the estimator $\widehat{H}_n(X)$ of $H$ by
    \begin{equation}
    \label{eq:H_est}
        \widehat{H}_n(X) := \frac{1}{2} - \frac{\log \left(V_{n,[0,T]}^2(X)\right) - \log \left(V_{\frac{n}{2},[0,T]}^2(X)\right)}{2 \log (2)}, \quad n\in2\N.
    \end{equation}
Estimator $\widehat{H}_n(X)$ is (weakly) consistent and satisfies the following non-central limit theorem:
    \[
    \left(\frac{1}{\delta_n} \right)^{2-2H'} (\widehat{H}_n(X) - H) \quad \xrightarrow[\substack{n\in 2 \N \\ n\to\infty}]{\P} \quad  \frac{d_{H',q}(2^{2-2H'}-1)}{2 \log(2) \int_{0}^{T} f_s^2 \d s}\int_{0}^{T} f_s^2 \d Z_s^{2H'-1,2}
    \]
\end{corollary}

\begin{proof}
    Using \autoref{thm: 2-vars distrib limit} and \autoref{lem:delta_method_for_p} yields
    \begin{multline}
        \label{eq: QV log}
        	 \left(\frac{1}{\delta_n} \right)^{2-2H'} \left(\int_{0}^{T} f_s^2 \d s \right)  \left[(2H-1) \log\left(\frac{1}{\delta_n}\right) + \log\left(V_{n,[0,T]}^2(X)\right) - \log\left(\int_{0}^{T} f_s^2 \d s \right) \right] \\
         \xrightarrow[n \to \infty]{\P} \quad  d_{H',q} \int_0^T f_s^2 \d Z_s^{2H'-1,2},
    \end{multline}
    and, consequently
    \begin{multline*}
         \left(\frac{1}{2\delta_n} \right)^{2-2H'} \left(\int_{0}^{T} f_s^2 \d s \right)  \left[(2H-1) \log\left(\frac{1}{2 \delta_n}\right)  + \log\left(V_{\frac{n}{2},[0,T]}^2(X)\right) - \log\left(\int_{0}^{T} f_s^2 \d s \right) \right] \\
         \qquad \xrightarrow[\substack{n\in 2\N\\ n \to \infty}]{\P} \quad  d_{H',q} \int_0^T f_s^2 \d Z_s^{2H'-1,2}.
    \end{multline*}
    Subtracting the two terms above concludes the proof.
\end{proof}

\begin{remark}
The estimator $\widehat{H}_n(X)$ of $H$ defined by \eqref{eq: QV log} is robust and fairly simple to use. From a practical point of view, its advantage is also that it does not explicitly depend on the Hermite order $q$ or on the noise intensity function $f$.
\end{remark}

\begin{remark}
 If we consider the scaled Rosenblatt process, i.e. $q=2$ and $f_s = \sigma$, $s\in [0,1]$, for some $\sigma>0$, we obtain from \autoref{cor: H estimator} the convergence
 \[
    \left(\frac{1}{\delta_n} \right)^{2-2H'} (\widehat{H}_n(X) - H) \quad \xrightarrow[\substack{n\in 2 \N \\ n\to\infty}]{\P} \quad \frac{2^\frac{3}{2}(2H-1)^\frac{1}{2}(2^{1-H}-1)}{H^\frac{1}{2}(H+1)\log(2)} Z_1^{H,2}, 
 \]
which corresponds to the asymptotic result for the estimator of $H$ in \cite{CouKriMas25}, except for the type of the convergence (in \cite{CouKriMas25}, the almost sure convergence was proved for the Rosenblatt process while here, we only obtain convergence in probability for a general Hermite process).
\end{remark}

\subsection{Estimation of the noise intensity}

Let us estimate the $L^2$-average of the noise intensity.

\begin{corollary}
\label{cor: avg intensity estim}
Assume that $q\geq 2$ or, if $q=1$, that $H>3/4$. Let also $f$ be a non-degenerate ($f\neq 0$), bounded, piecewise $\alpha$-H\"older continuous function on $[0,T]$ for some $\alpha>2-2H'$ where $H'$ is defined by \eqref{eq:H'}. Denote the $L^2$-average of the noise intensity by
    \[
    	\sigma_2 := \left(\int_{0}^{T} f_s^2 \d s\right)^\frac{1}{2} 
    \]
and define its estimator $\widehat{\sigma}_{2,n}(X)$ by
    \begin{equation}
        \label{eq: sigma est def}
        \widehat{\sigma}_{2,n} := \exp\left[ \left(\widehat{H}_n(X) - \frac{1}{2}\right) \log\left(\frac{1}{\delta_n}\right) + \frac{1}{2}  \log \left(V_{n,[0,T]}^2(X)\right) \right]
    \end{equation}
where $\widehat{H}_n(X)$ is the estimator of the Hurst index defined by \eqref{eq:H_est}. Then there is the following joint convergence of the two estimators:
    \begin{equation}\label{eq: convergence H sigma}
       \left(\frac{1}{\delta_n} \right)^{2-2H'} \begin{bmatrix} 
             \widehat{H}_n(X) - H \\
            \frac{1}{\log(1/\delta_n)} (\widehat{\sigma}_{2,n}(X)  - \sigma_2) \\
        \end{bmatrix} 
        \quad \xrightarrow[\substack{n\in 2 \N \\ n\to\infty}]{\P} \quad
        \begin{bmatrix} 
            1\\
            \sigma_2\\
        \end{bmatrix}
         \frac{d_{H',q}(2^{2-2H'}-1)}{2 \log(2) \int_{0}^{T} f_s^2 \d s}\int_{0}^{T} f_s^2 \d Z_s^{2H'-1,2}.
    \end{equation}
\end{corollary}

\begin{proof}
Convergence \eqref{eq: convergence H sigma} of the random vectors follows if its individual components converge. The first component converges by \autoref{cor: H estimator}. To prove the convergence of the second component note that we have 
\[
        \left(\frac{1}{\delta_n} \right)^{2-2H'} \frac{1}{\log(1/\delta_n)} \left[(2\widehat{H}_n(X)-1) \log\left(\frac{1}{\delta_n}\right) + \log\left(V_{n,[0,T]}^2(X)\right) - \log\left(\int_{0}^{T} f_s^2 \d s \right) \right] = a_n + b_n
\]
where 
	\begin{align*}
		a_n & :=  2\left(\frac{1}{\delta_n} \right)^{2-2H'}  (\widehat{H}_n(X) - H), \\
		b_n & := \frac{1}{\log(1/\delta_n)} \left(\frac{1}{\delta_n} \right)^{2-2H'}  \left[(2H-1) \log\left(\frac{1}{\delta_n}\right) + \log\left(V_{n,[0,T]}^2(X)\right) - \log\left(\int_{0}^{T} f_s^2 \d s \right) \right].
	\end{align*}
We clearly have 
    \[
    	a_n \quad \xrightarrow[n\to\infty]{\P}\quad  \frac{d_{H',q}(2^{2-2H'}-1)}{\log(2) \int_{0}^{T} f_s^2 \d s}\int_{0}^{T} f_s^2 \d Z_s^{2H'-1,2}
     \] 
by \autoref{cor: H estimator} and we also have
	\[
	 	b_n = O_{\P}\left(\frac{1}{\log(1/\delta_n)}\right), \quad n\to\infty,
	\] 
by convergence \eqref{eq: QV log}. Consequently, there is the convergence
    \begin{multline}
    \label{eq: plugin convergence}
     \left(\frac{1}{\delta_n} \right)^{2-2H'} \frac{1}{\log(1/\delta_n)} \left[ \log\left(\frac{1}{\delta_n^{2\widehat{H}_n(X)-1}}V_{n,[0,T]}^2(X)\right) - \log\left(\int_{0}^{T} f_s^2 \d s \right) \right] \\
    \qquad \xrightarrow[\substack{n\in 2 \N \\ n\to\infty}]{\P} \quad  \frac{d_{H',q}(2^{2-2H'}-1)}{\log(2) \int_{0}^{T} f_s^2\d s}\int_{0}^{T} f_s^2 \d Z_s^{2H'-1,2}
    \end{multline}
and the claim now follows by the use of \autoref{lem:delta_method_for_p} with function $g(x) = \exp(x/2)$, $x\in\R$.
\end{proof}

\begin{remark}
Up to the convergence type, \autoref{cor: avg intensity estim} generalizes \cite[Theorem 3.2]{CouKriMas25} where only the scaled Rosenblatt process is considered. Here we consider not only general Hermite process but also non-constant noise intensity function.
\end{remark}

In what follows, we investigate non-parametric estimation of the noise intensity function.

\begin{corollary}
\label{thm: weak consistency of noise estimator}
Assume that $q\geq 2$ or, if $q=1$, that $H>3/4$. Let also $f$ be a non-degenerate $(f\neq 0)$, bounded, piecewise $\alpha$-H\"older continuous function on $[0,T]$ for some $\alpha>2-2H'$ where $H'$ is defined by \eqref{eq:H'}. Then the estimator $\widehat{f}_n(X)$ defined by 
	\begin{equation}
    \label{eq: estimator}
    \widehat{f}_n(X)(t):=\sum_{j=0}^{K_n-1} \left(\frac{\cfrac{1}{\delta_n^{2 \widehat{H}_n(X)-1}}V_{n,[T_j^n,T_{j+1}^n]}^2\left( X \right)}{T_{j+1}^n-T_j^n}\right)^{1/2} \bm{1}_{[T_j^n,T_{j+1}^n)}(t), \quad t\in [0,T], n\in\N,
\end{equation}
where $\widehat{H}_n(X)$ is the estimator of the Hurst parameter $H$ defined by \eqref{eq:H_est}, is a (weakly) consistent estimator of the noise intensity function $|f|$ with respect to the the $L^{2}$-distance, i.e.\ there is the following convergence:
	\begin{equation}
	 \|\widehat{f}_n(X)-|f|\|_{L^2(0,T)} \xrightarrow[n \to \infty]{\P} 0.
	\end{equation}
\end{corollary}

\begin{proof}
    The convexity of the function $f(x) = x^2$, $x\in\R$, implies  the inequality 
    \[
    	 |x - y|^{2} \leq   |x^2-y^2|, \quad x,y\geq 0. 
   	\]
    By using this fact, we obtain   
    \begin{align*}
         \int_0^T\left|\widehat{f}_n(X)(t) - |f_t|\right|^{2}\d{t} &  \\
         & \hspace{-3cm} \leq  \int_0^T\left|\widehat{f}_n(X)(t)^{2} - f_t^{2}\right|\d{t}\\
        &  \hspace{-3cm} = \sum_{j=0}^{K_n-1} \int_{T_j^n}^{T_{j+1}^n}\left|\frac{\cfrac{1}{\delta_n^{2 \widehat{H}_n(X)-1}}V_{n,[T_j^n,T_{j+1}^n]}^2\left( X \right)}{T_{j+1}^n-T_j^n}- f_t^{2}\right|\d{t}\\
        & \hspace{-3cm} =\sum_{j=0}^{K_n-1} \int_{T_j^n}^{T_{j+1}^n}\left| \frac{\cfrac{1}{\delta_n^{2 \widehat{H}_n(X)-1}}V_{n,[T_j^n,T_{j+1}^n]}^2\left( X \right)}{T_{j+1}^n-T_j^n} - \frac{\int_{T_j^n}^{T_{j+1}^n}f_s^2\d{s}}{T_{j+1}^n-T_j^n}   +  \frac{\int_{T_j^n}^{T_{j+1}^n}f_s^{2}\d{s}}{T_{j+1}^n-T_j^n}  - f_t^{2}\right|\d{t}\\
        & \hspace{-3cm} \leq a_n + b_n
    \end{align*}
where 
	\begin{align*}
		a_n & := \sum_{j=0}^{K_n-1} \left| \cfrac{1}{\delta_n^{2 \widehat{H}_n(X)-1}}V_{n,[T_j^n,T_{j+1}^n]}^2\left( X \right) - \int_{T_j^n}^{T_{j+1}^n}f_s^{2}\d{s} \right|, \\
		b_n & := \sum_{j=0}^{K_n-1}\int_{T_j^n}^{T_{j+1}^n}\left| \frac{\int_{T_j^n}^{T_{j+1}^n}f_s^{2}\d{s}}{T_{j+1}^n-T_j^n}  - f_t^2\right|\d{t}.
	\end{align*}
    By \autoref{thm: vars joint convergence} and \autoref{cor: H estimator}, we immediately obtain the convergence
    \[ a_n \xrightarrow[n \to \infty]{\P} 0,\]
    while the convergence
    \[ b_n \xrightarrow[n \to \infty]{} 0 \]
    follows from piecewise uniform continuity of $f^2$.
\end{proof}

\begin{remark}
We note that in order to use the estimators from \autoref{cor: avg intensity estim} and \autoref{thm: weak consistency of noise estimator}, one has to compute the estimate of the Hurst parameter $\widehat{H}_n(X)$ first. Then, however, the estimators $\widehat{\sigma}_{2,n}(X)$ and $\widehat{f}_{n}(X)$ can be immediately computed as the Hermite order $q$ is not needed.
\end{remark}

\subsection{Estimation of Hermite order $q$}

In this section, we let $f\equiv \sigma$ where $\sigma\in (0,\infty)$ so that $X=\sigma Z^{H,q}$. We assume that $\sigma$ is known and our aim is to estimate the Hermite order $q$. 

\begin{corollary}
\label{thm: est H' and q} 
Assume that $q\geq 2$ or, if $q=1$, that $H>3/4$. Denote
    \[
    \widetilde{F}_n(X) := (2\widehat{H}_n(X)-1) \log\left(1/\delta_n\right) + \log(V_{n,[0,T]}^2(X)) - \log(T\sigma^2 ), \quad n\in 2\N,
    \]
where $\widehat{H}_n(X)$ is the estimator of the Hurst index $H$ defined by \eqref{eq:H_est} and define the estimator $\widehat{H}'_n(X)$ of $H'$ by
    \[
    \widehat{H}'_n(X):=1+ \frac{\log(\log(1/(2\delta_n))) - \log(\log(1/\delta_n)) +\log (|\widetilde{F}_n(X)|) - \log (|\widetilde{F}_{n/2}(X)|)}{2 \log(2)}, \quad n\in 4\N.
    \]
    Then there is the convergence
    \begin{equation*}
     \widehat{H}'_n(X)\xrightarrow[\substack{ n\in 4\N\\ n \to \infty}]{\P} H'
    \end{equation*}
    and, consequently,
    \begin{equation*}
     \hat{q}_n(X) := \frac{\widehat{H}_n(X) -1}{\widehat{H}'_n(X) - 1} \quad\xrightarrow[\substack{n\in 4\N\\ n \to \infty}]{\P}\quad \frac{H-1}{H'-1} = q.
    \end{equation*}
\end{corollary}

\begin{proof}
    Denote
    \[
    \widetilde{Y}:=\frac{d_{H',q}(2^{2-2H'}-1)}{ T\log(2)}Z_T^{2H'-1,2}.
    \]
    By \eqref{eq: plugin convergence}, we have the convergence
    \[
    \left(\frac{1}{\delta_n} \right)^{2-2H'} \frac{1}{\log(1/\delta_n)} \, \tilde{F}_n(X)  \xrightarrow[\substack{n\in 4\N\\ n \to \infty}]{\P} \tilde{Y}.
    \]
    By the continuous mapping theorem we obtain the following convergences:
    \begin{align*}
        & (2-2H') \log\left(1/ \delta_n \right) - \log(\log(1/\delta_n)) +\log (|\tilde{F}_n(X)|) \xrightarrow[\substack{n\in 2\N\\ n \to \infty}]{\P} \log(|\tilde{Y}|), \\
        & (2-2H') \left(\log\left(1/ \delta_n\right) - \log (2)\right) - \log(\log(1/(2\delta_n))) +\log (|\tilde{F}_{n/2}(X)|) \xrightarrow[\substack{n\in 4\N\\ n \to \infty}]{\P} \log(|\tilde{Y}|).
    \end{align*}
    Subtracting the two quantities yields
    \[
    (2-2H') \log 2 + \log(\log(1/(2\delta_n))) - \log(\log(1/\delta_n)) +\log (|\tilde{F}_n(X)|) - \log (|\tilde{F}_{n/2}(X)|)  \xrightarrow[\substack{n\in 4\N\\ n \to \infty}]{\P} 0, 
    \]
    and, after rearranging,
    \[
    1+ \frac{\log(\log(1/(2\delta_n))) - \log(\log(1/\delta_n)) +\log (|\tilde{F}_n(X)|) - \log (|\tilde{F}_{n/2}(X)|)}{2 \log 2} \xrightarrow[\substack{n\in 4\N\\ n \to \infty}]{\P} H'.
    \]
\end{proof}

\begin{remark}
Note that the estimator $\widehat{q}_n(X)$ of order $q$ can be evaluated only if $\sigma$ is known. The situation where $\sigma$ is unknown is discussed in \autoref{app:unknown_scaling}.
\end{remark}

\section{Estimation for non-linear SDEs with additive Hermite noise}
\label{sec:SDEs}

In what follows, we aim to generalize the results of the previous sections to solutions of stochastic differential equations (SDEs). Consider the SDE
\begin{equation}
	\label{eq:non-lin_SDE}
	 Y_t = Y_0 +  \int_0^t g(Y_s)\d{s} + \int_0^t f_s \d Z_s^{H,k}, \qquad t\in [0,T].
\end{equation}
Here, we assume that $Z^{H,q}$ is a $q$\textsuperscript{th} order Hermite process with Hurst parameter $H\in (1/2,1)$, $Y_0\in\R$, $f:[0,T] \to \R$ is a deterministic, Borel measurable, and bounded function, and $g:\R \to \R$ is a deterministic, Borel measurable function that satisfies the following two conditions:
\begin{enumerate}[label=(g\arabic*)]
	\item\label{f1} Function $g$ is locally Lipschitz, i.e.\ for every $N\in\N$ there exists a constant $K_N\in (0,\infty)$ such that for every $x,y\in\R$, $|x|+|y|\leq N$, it holds that
    \[ |g(x)-g(y)|\leq K_N |x-y|;\]
	\item\label{f2} There exists a (Lyapunov) function $V\in C^1(\R)$ such that
		\begin{equation*}
		\label{f2a} \lim_{R\to\infty} \inf_{|x|>R} V(x)=\infty
		\end{equation*}
 and for which there is a constant $K\in (0,\infty)$ and a continuous function $h:[0,\infty)\to [0,\infty)$ such that the following inequality holds for every $x,z\in\R$:
		\begin{equation*}
		\label{f2b}
		 V'(x)g(x+z) \leq K V(x) + h(|z|).
		\end{equation*}
\end{enumerate}

By repeating the approach from \cite{CouKriMas25}, that is based on the connection of solutions to SDEs with additive noise with pathwise solutions to the corresponding random differential equations, we obtain the following result: 

\begin{theorem}
There exists a unique solution to equation \eqref{eq:non-lin_SDE} with almost every sample path in $C^{\nu}([0,T])$ for every $\nu\in (0,H)$.
\end{theorem} 

The following two results show that the behavior of the quadratic variation of the solution $Y$ to equation \eqref{eq:non-lin_SDE} corresponds, under some mild additional assumptions, to the behavior of the quadratic variation of the driving Hermite process.

\begin{theorem}
\label{thm: vars joint convergence SDE}
Let $Y$ be the solution to equation \eqref{eq:non-lin_SDE}. Then
\[
	\sum_{j=0}^{K_n-1} \left| \frac{1}{\delta_n^{2H-1}} V_{n,[T_j^n, T_{j+1}^n]}^2(Y) - \int_{T_j^n}^{T_{j+1}^n} f_s^2\d{s}\right| \quad\xrightarrow[n\to\infty]{\P} \quad 0.
\]
\end{theorem}

\begin{proof}
Let us decompose the solution $Y$ to the (smooth) drift and the (rough) stochastic part, i.e.\ define processes $Y^S=(Y_t^S)_{t\in [0,T]}$ and $Y^R=(Y_t^R)_{t\in [0,T]}$ by
	\[ Y_t^S := \int_0^t g(Y_s)\d{s} \quad\mbox{and}\quad Y_t^R := \int_0^t f_s\d{Z}_{s}^{H,q}\]
for $t\in [0,T]$ so that $Y=Y_0+ Y^S + Y^R$. It follows by \autoref{thm: vars joint convergence} that it suffices to prove that
	\[ 
		\delta_n^{1-2H}\sum_{j=0}^{K_n-1} \left|V_{n,[T_j^n,T_{j+1}^n]}^2 (Y)-V_{n,[T_j^n,T_{j+1}^n]}^2 (Y^R)\right| \xrightarrow[n\to\infty]{\P} 0.
	\]
To this end, let $2H-1<\nu<H$ (such choice of $\nu$ is possible as the equality $2H-1<H$ is equivalent to $H<1$ which is satisfied). Consider the $\nu$-H\"older continuous version of $Y^R$ (denoted again by $Y^R$). We have that 
	\begin{align*}
	  \left| V_{n,[T_j^n,T_{j+1}^n]}^2(Y) - V_{n,[T_j^n,T_{j+1}^n]}^2(Y^R)\right| & \leq \sum_{i: t_i^n\in [T_j^n, T_{j+1}^n)} \left|\left(Y_{t_i^n,t_{i+1}^n}^S + Y_{t_i^n,t_{i+1}^n}^R\right)^2 - \left(Y_{t_i^n,t_{i+1}^n}^R\right)^2  \right| \\
	  	& \leq  \sum_{i: t_i^n\in [T_j^n, T_{j+1}^n)} \left(Y_{t_i^n,t_{i+1}^n}^S\right)^2 + 2\sum_{i: t_i^n\in [T_j^n, T_{j+1}^n)}\left|Y_{t_i^n,t_{i+1}^n}^S\right|\left|Y_{t_i^n,t_{i+1}^n}^R\right| \\
	  	& \leq \|g(Y)\|_\infty^2 D_n \delta_n^2 + 2\|g(Y)\|_\infty \|Y^R\|_{C^\nu} D_n\delta_n^{1+\nu}
	\end{align*} 
holds for $n\in\N$ and $j\in\{0,1,\ldots, K_n-1\}$ almost surely so that 	
	\[
		\delta_n^{1-2H}\sum_{j=0}^{K_n-1} \left|V_{n,[T_j^n,T_{j+1}^n]}^2 (Y)-V_{n,[T_j^n,T_{j+1}^n]}^2 (Y^R)\right| \lesssim K_nD_n \delta_n^{2-2H+\nu} = T \delta_n^{1-2H+\nu} \xrightarrow[n\to\infty]{a.s.} 0.
	\]
\end{proof}

\begin{theorem}
\label{thm: 2-vars distrib limit SDE}
Assume that $q\geq 2$ and assume that $f$ is additionally bounded, piecewise $\alpha$-H\"older continuous function on $[0,T]$ for some $\alpha>2-2H'$. Let $Y$ be the solution to equation \eqref{eq:non-lin_SDE}. Then 
	\[
		\frac{1}{\delta_n^{2-2H'}} \left( \frac{1}{\delta_n^{2H-1}} V_{n,[0,T]}^2(Y) - \int_0^T f_s^2\d{s}\right) \quad \xrightarrow[n\to\infty]{\P} \quad d_{H',q} \int_0^T f_s^2\d{Z}_s^{2H'-1,2}.
	\]
\end{theorem}

\begin{proof}
Let us decompose the solution $Y$ to the (smooth) drift $Y^S$ and the (rough) stochastic part $Y^R$ as in the proof of \autoref{thm: vars joint convergence SDE}. It follows by \autoref{thm: 2-vars distrib limit} that it suffices to prove that
	\[
		\delta_n^{2H'-2H-1} \left| V_{n,[0,T]}^2(Y) - V_{n,[0,T]}^2(Y^R)\right| \xrightarrow[n\to\infty]{\P} 0.
	\]  
To this end, note that the estimate
	\[
	  \left| V_{n,[0,T]}^2(Y) - V_{n,[0,T]}^2(Y^R)\right| \leq  \sum_{i=0}^{n-1} \left(Y_{t_i^n,t_{i+1}^n}^S\right)^2 + 2\left|\sum_{i=0}^{n-1}Y_{t_i^n,t_{i+1}^n}^SY_{t_i^n,t_{i+1}^n}^R\right| 
	\]
holds for every $n\in\N$ almost surely. For the first term, we have that the estimate
	\[ 
			\sum_{i=0}^{n-1} \left( Y_{t_i^n, t_{i+1}^n}^S\right)^2 \leq \|g(Y)\|_{\infty}^2 T\delta_n 
	\] 
holds for every $n\in\N$ almost surely so that 
	\[
		\delta_n^{2H'-2H-1}\sum_{i=0}^{n-1} \left( Y_{t_i^n, t_{i+1}^n}^S\right)^2 \xrightarrow[n\to\infty]{a.s.}0.
	\]	
holds as $H'>H$ by the assumption that $q\geq 2$. For the second term, let $\nu\in (1/2,H)$ and let $\omega\in\Omega$ be such that both $Y(\omega)$ and $Y^R(\omega)$ are $\nu$-H\"older continuous functions. By the mean value theorem, we have the equality 
	\[
		\sum_{i=0}^{n-1} Y_{t_i^n,t_{i+1}^n}^S(\omega) Y_{t_i^n,t_{i+1}^n}^R(\omega) = \delta_n \sum_{i=0}^{n-1} g(Y_{s_i^n}(\omega))Y_{t_i^n, t_{i+1}^n}^R(\omega)
	\]
for some $s_i^n\in [t_i^n, t_{i+1}^n]$. Note that there is $M(\omega)\in\N$ such that
	\[
	 	|Y_s(\omega)| + |Y_t(\omega)| \leq 2\|Y(\omega)\|_{C^\nu([0,T])} \leq M(\omega), \quad s,t\in [0,T],
	\]
and thus, by assumption \ref{f1}, we have the estimate 
	\[ 
		|g(Y_s(\omega)) - g(Y_t(\omega))| \leq K_{M(\omega)} |Y_s(\omega)-Y_t(\omega)| \leq K_{M(\omega)} \|Y(\omega)\|_{C^\nu([0,T])} |s-t|^\nu, \quad s,t\in [0,T],
	\]
which proves that also $g(Y(\omega))\in C^\nu([0,T])$. As $2\nu>1$, it follows from \cite{You36} that there is the convergence
	\[ 
		\sum_{i=0}^{n-1} g(Y_{s_i^n}(\omega))Y_{t_i^n, t_{i+1}^n}^R(\omega) \quad \xrightarrow[n\to\infty]{} \quad\int_0^T g(Y_s(\omega))\circ \d{Y^R_s(\omega)}
	\]
where the integral on the right-hand side is the Young integral. Consequently, we have that 
	\[ 
		\delta_n^{2H'-2H-1} \left|\sum_{i=0}^{n-1} Y_{t_i^n, t_{i+1}^n}^SY_{t_i^n, t_{i+1}^n}^R\right| \xrightarrow[n\to\infty]{a.s.} 0.
	\]
\end{proof}

In what follows, we focus on the applications of \autoref{thm: vars joint convergence SDE} and \autoref{thm: 2-vars distrib limit SDE} to (parametric and non-parametric) estimation for the SDE \eqref{eq:non-lin_SDE}. We assume that we have a single trajectory of the solution $Y$ to \eqref{eq:non-lin_SDE} on interval $[0,T]$ observed discretely at points $\{t_i^n\}_{i=0}^n$. It turns out that, in the non-Gaussian case ($q\geq 2$), the results from section \ref{sec:estimation_wiener_integral} transfer to this situation.

We begin with a (weakly) consistent estimator of the Hurst parameter of the driving Hermite process. As before, the estimator does not depend on the other parameters (i.e.\ the drift function $g$, the diffusion coefficient $f$, and the order of the Hermite process $q$) and it is easily implemented.

\begin{corollary}
\label{cor: H estimator SDE}
Assume that $q\geq 2$ and assume that $f$ is additionally non-degenerate ($f\neq 0$), bounded, piecewise $\alpha$-H\"older continuous function on $[0,T]$ for some $\alpha>2-2H'$. Let $Y$ be the solution to equation \eqref{eq:non-lin_SDE}. Define the estimator $\widehat{H}_n$ of $H$ by 
	\begin{equation}
	\label{eq:est_H_SDE}
		\widehat{H}_n(Y) := \frac{1}{2} - \frac{\log\left(V_{n,[0,T]}^2(Y)\right)-\log \left( V_{\frac{n}{2},[0,T]}^2(Y)\right)}{2\log (2)}, \quad n\in 2\N.
	\end{equation}
Estimator $\widehat{H}_n(Y)$ is (weakly) consistent and satisfies the following non-central limit theorem:
	\[ 
		\left(\frac{1}{\delta_n}\right)^{2-2H'} (\widehat{H}_n(Y)-H) \quad\xrightarrow[n\to\infty]{\P} \quad \frac{d_{H',q}(2^{2-2H'}-1)}{2\log (2) \int_0^T f_s^2\d{s}}\int_0^T f_s^2 \d{Z}_s^{2H'-1,2}.
	\] 
\end{corollary}

\begin{proof}
\autoref{cor: H estimator SDE} follows from \autoref{thm: 2-vars distrib limit SDE} in the same manner as \autoref{cor: H estimator} follows from \autoref{thm: 2-vars distrib limit}.
\end{proof}

Let us now continue with a (parametric) estimator of the $L^2$-average of the noise intensity $|f|$ and with a (non-parametric) estimator of the noise intensity $|f|$. Again, these estimators do not depend on the drift function $g$ or the parameters $H$ and $q$ but they require that the estimator $\widehat{H}_n(Y)$ is computed first. 

\begin{corollary}
\label{cor:est_s2_SDE}
Assume that $q\geq 2$ and that $f$ is additionally non-degenerate ($f\neq 0$), bounded, piecewise $\alpha$-H\"older continuous function on $[0,T]$ for some $\alpha>2-2H'$. Let $Y$ be the solution to equation \eqref{eq:non-lin_SDE}. Denote the $(L^2)$-average of the noise intensity by 
	\[ 
		\sigma_2 := \left(\int_0^T f_s^2\d{s}\right)^\frac{1}{2}
	\]
and define its plug-in estimator 
	\[ 
		\widehat{\sigma}_{2,n}(Y) := \exp\left[ \left(\widehat{H}_n(Y)-\frac{1}{2}\right) \log\left(\frac{1}{\delta_n}\right) + \frac{1}{2} \log \left( V_{n,[0,T]}^2(Y)\right)\right],
	\]
where $\widehat{H}_n(Y)$ is the estimator of Hurst index defined by \eqref{eq:est_H_SDE}. Then there is the following joint convergence of the two estimators: 
	\[
		\left(\frac{1}{\delta_n}\right)^{2-2H'}
		\begin{bmatrix}
		  \widehat{H}_n(Y)-H \\
		 	\frac{1}{\log(1/\delta_n)} (\widehat{\sigma}_{2,n}(Y) - \sigma_2) 
		\end{bmatrix}\quad \xrightarrow[n\to\infty]{\P} \quad \frac{d_{H',q}(2^{2-2H'}-1)}{2\log(2) \int_0^T f_s^2\d{s}} \int_0^T f_s^2\d{Z}_s^{2H'-1,2}.
	\]
\end{corollary}

\begin{proof}
\autoref{cor:est_s2_SDE} follows from \autoref{cor: H estimator SDE} in the same manner as \autoref{cor: avg intensity estim} follows from \autoref{cor: H estimator}.
\end{proof}

\begin{corollary}
\label{cor: weak consistency of noise estimator SDE} 
Assume that $q\geq 2$ and assume that $f$ is additionally non-degenerate ($f\neq 0$), bounded, piecewise $\alpha$-H\"older continuous function on $[0,T]$ for some $\alpha>2-2H'$. Let $Y$ be the solution to equation \eqref{eq:non-lin_SDE}. Then the estimator $\widehat{f}_n$ defined by 
	\begin{equation*}
    \widehat{f}_n(Y)(t):=\sum_{j=0}^{K_n-1} \left(\frac{\frac{1}{\delta_n^{2 \widehat{H}_n(Y)-1}}V_{n,[T_j^n,T_{j+1}^n]}^2\left( Y \right)}{T_{j+1}^n-T_j^n}\right)^{1/2} \bm{1}_{[T_j^n,T_{j+1}^n)}(t), \quad t\in [0,T], n\in\N,
\end{equation*}
where $\widehat{H}_n(Y)$ is the estimator of the Hurst index $H$ defined by \eqref{eq:est_H_SDE} is a (weakly) consistent estimator of the noise intensity function $|f|$ with respect to the the $L^{2}$-distance, i.e.\ there is the following convergence:
	\begin{equation}
	 \|\widehat{f}_n(Y)-|f|\|_{L^2(0,T)} \xrightarrow[n \to \infty]{\P} 0.
	\end{equation}
\end{corollary}

\begin{proof}
\autoref{cor: weak consistency of noise estimator SDE} follows from \autoref{thm: vars joint convergence SDE} and \autoref{cor: H estimator SDE} in the same manner as \autoref{thm: weak consistency of noise estimator} follows from \autoref{thm: vars joint convergence} and \autoref{cor: H estimator}.
\end{proof}

We now turn our attention to the estimation of parameter $q$.

\begin{corollary}
\label{thm: est H' and q SDE} 
Assume that $q\geq 2$ and let $Y$ be the solution to equation \eqref{eq:non-lin_SDE} with $f \equiv \sigma$ where $\sigma\in (0,\infty)$. Denote
    \[
    \widetilde{F}_n(Y) := (2\widehat{H}_n(Y)-1) \log\left(1/\delta_n\right) + \log(V_{n,[0,T]}^2(Y)) - \log(T\sigma^2 ), \quad n\in 2\N,
    \]
where $\widehat{H}_n(Y)$ is the estimator of the Hurst index $H$ defined by \eqref{eq:est_H_SDE} and define the estimator $\widehat{H}'_n(Y)$ of $H'$ by
    \[
    \widehat{H}'_n(Y):=1+ \frac{\log(\log(1/(2\delta_n))) - \log(\log(1/\delta_n)) +\log (|\widetilde{F}_n(Y)|) - \log (|\widetilde{F}_{n/2}(Y)|)}{2 \log(2)}, \quad n\in 4\N.
    \]
    Then there is the convergence
    \begin{equation*}
     \widehat{H}'_n(Y)\xrightarrow[\substack{ n\in 4\N\\ n \to \infty}]{\P} H'
    \end{equation*}
    and, consequently,
    \begin{equation*}
     \hat{q}_n(Y) := \frac{\widehat{H}_n(Y) -1}{\widehat{H}'_n(Y) - 1} \quad\xrightarrow[\substack{n\in 4\N\\ n \to \infty}]{\P}\quad \frac{H-1}{H'-1} = q.
    \end{equation*}
\end{corollary}

\begin{proof}
\autoref{thm: est H' and q SDE} follows from \autoref{cor: H estimator SDE} in a similar manner as \autoref{thm: est H' and q} follows from \autoref{cor: H estimator}.
\end{proof}

\bibliographystyle{siam}
\bibliography{bibi,bibi_loc}

\appendix

\section{Estimation of $q$ when parameter $\sigma$ is unknown}
\label{app:unknown_scaling}

In this appendix, we consider the situation in which we have a single trajectory of $X=\sigma Z^{H,q}$ observed discretely at points $\{t_i^n\}_{i=0}^n$ but none of the parameters are known. We wish to estimate the Hermite order $q$. To this end, we first estimate $\sigma$ from all data on $[0,T]$ and then insert the estimate into the estimator of $q$ from \autoref{thm: est H' and q} evaluated separately on the sub-intervals $[T_j^n,T_{j+1}^n]$, $j\in \{0,1,\ldots, K_{n-1}\}$. More precisely, there is the following result.

\begin{theorem}
    \label{thm: est H' sigma and q} 
    Let $\sigma\in (0,\infty)$, and let $X=\sigma Z^{H,q}$. Assume that $q\geq 2$ or, if $q=1$, that $H>3/4$. Assume also that\footnote{For example, we can choose $\Delta_n = \delta_n^\alpha$ for some $\alpha\in (0,1)$.}
    \begin{equation}
        \label{eq: assum Delta delta}
        \Delta_n^{2-2H'} \log(1/\delta_n) \xrightarrow[n \to \infty]{} 0.
    \end{equation}
    Denote
     \begin{equation}
     \label{eq:plugin_F}
     	\widehat{F}_{n,[T_{j_n}^n,T_{j_n+1}^n]}(X):= (2\widehat{H}_n(X)-1) \log(1/\delta_n) + \log(V_{n,[T_{j_n}^n,T_{j_n+1}^n]}^2(X)) - \log(\Delta_n\widehat{\sigma}_{2,n}(X)^2)
   \end{equation}
 and define the estimator $\widehat{H}'_{n,j_n}(X)$ of $H'$ by
    \[
    \widehat{H}_{n,j_n}'(X):= 1+ \frac{ \log (|\widehat{F}_{n,[T_{j_n}^n,T_{j_n+1}^n]}(X)|) - \log (|\widehat{F}_{n/2,[T_{j_n}^n,T_{j_n+1}^n]}(X)|)}{2 \log(2)}
    \]
and the estimator $\widehat{q}_{n,j_n}(X)$ of $q = \frac{H-1}{H'-1}$ by 
	\[ 
		\widehat{q}_{n,j_n}(X) := \frac{\widehat{H}_n(X) -1}{\widehat{H}_{n,j_n}'(X)  - 1}
	\]
for $j_n\in\{0,1, \ldots, K_n-1\}$ and $n\in2\N$. 
    Then there is the convergence
    \begin{equation*}
     \widehat{H}_{n,j_n}'(X) \xrightarrow[\substack{n\in 2\N\\ n \to \infty}]{\P} H',
    \end{equation*}
    and, consequently, also the convergence
    \begin{equation*}
     \widehat{q}_{n,j_n}(X)\xrightarrow[\substack{n\in 2\N\\ n \to \infty}]{\P} q.
    \end{equation*}    
\end{theorem}

\begin{proof}
	Let us denote 
   	\[
   	F_{n,[T_{j_n}^n,T_{j_n+1}^n]}(X):= (2H-1) \log\left(1/\delta_n\right) + \log\left(V_{n,[T_{j_n}^n,T_{j_n+1}^n]}^2(X)\right) - \log\left(\sigma^2 (T_{j_n+1}^n - T_{j_n}^n)\right)
    \]
   for $j_n\in \{0,1,\ldots, K_n-1\}$ and $n\in2 \N$. Self-similarity and stationarity of the increments of $X = \sigma Z^{H,q}$ yields the equality in distribution
    \begin{align*}
        &\begin{pmatrix}
        F_{n,[T_{j_n}^n,T_{j_n+1}^n]}(X)\\
        F_{\frac{n}{2},[T_{j_n}^n,T_{j_n+1}^n]}(X) 
        \end{pmatrix}
        \eqd \begin{pmatrix}
        F_{\frac{\Delta_n}{\delta_n},[0,T]}(X)\\
        F_{\frac{\Delta_n}{2\delta_n},[0,T]}(X) 
        \end{pmatrix}
    \end{align*}
for any $j_n\in \{0,1,\ldots, K_n-1\}$ and $n\in2\N$. By \eqref{eq: QV log}, we have the convergence
    \begin{align*}
        &\begin{pmatrix}
        \left(\frac{\Delta_n}{\delta_n} \right)^{2-2H'} F_{\frac{\Delta_n}{\delta_n},[0,T]}(X)\\
        \left(\frac{\Delta_n}{2\delta_n} \right)^{2-2H'}F_{\frac{\Delta_n}{2\delta_n},[0,T]}(X)
        \end{pmatrix}
        \xrightarrow[\substack{n\in 2\N\\ n \to \infty}]{\P}
        \begin{pmatrix}
        Y\\
        Y
        \end{pmatrix},
    \end{align*}
where 
    \[
    Y :=  \frac{d_{H',q}}{T} Z_T^{2H'-1,2}.
    \]
   Further note that it follows from \autoref{cor: avg intensity estim} that 
    \[\left(\frac{\Delta_n}{\delta_n} \right)^{2-2H'} \left( F_{n,[T_{j_n}^n,T_{j_n+1}^n]}(X) - \widehat{F}_{n,[T_{j_n}^n,T_{j_n+1}^n]}(X)\right)  = O_{\P}(\Delta_n^{2-2H'} \log(1/\delta_n))
    \]
for any $j_n\in\{0,1,\ldots, K_n-1\}$, $n\in 2\N$. In particular, a Slutsky-type argument provides us with the following joint convergence in distribution:
    \begin{align*}
        &\begin{pmatrix}
        \left(\frac{\Delta_n}{\delta_n} \right)^{2-2H'} \widehat{F}_{n,[T_{j_n}^n,T_{j_n+1}^n]}(X)\\
        \left(\frac{\Delta_n}{2\delta_n} \right)^{2-2H'} \widehat{F}_{\frac{n}{2},[T_{j_n}^n,T_{j_n+1}^n]}(X)
        \end{pmatrix}
        \xrightarrow[\substack{n\in 2\N\\n \to \infty}]{\D}
        \begin{pmatrix}
        Y\\
        Y
        \end{pmatrix}.
 	\end{align*}
    The rest of the proof is analogous to the proof of \autoref{thm: est H' and q} with the convergence in probability therein replaced by the convergence in distribution. However, since the resulting limit is a deterministic constant, i.e.\
    \begin{equation}
        1+ \frac{ \log |\widehat{F}_{n,[T_{j_n}^n,T_{j_n+1}^n]}(X)| - \log |\widehat{F}_{\frac{n}{2},[T_{j_n}^n,T_{j_n+1}^n]}(X)|}{2 \log(2)} \xrightarrow[\substack{n\in2\N\\ n \to \infty}]{\D} H',
    \end{equation}
    we automatically get the convergence in probability to the same constant $H'$.
\end{proof}

\begin{remark}
Comparing estimator $\widehat{H}'_n$ from \autoref{thm: est H' and q} to estimator $\widehat{H}_{n,j_n}'$ from \autoref{thm: est H' sigma and q}, we note that there is the additional term
\[	
	\log(\log(1/(2\delta_n))) - \log(\log(1/\delta_n))
\] 
which, however, converges to zero as $n\to \infty$.
\end{remark}

\begin{remark}
It is not quite clear, how to combine the estimators $\widehat{q}_{n,j_n}(X)$ of $q$ evaluated separately on the sub-intervals $[T_j^n, T_{j+1}^n]$ into a single estimator. A natural choice would be to take the average, however, it is not at all clear whether the estimator obtained in this manner would converge to the true value of $q$.
\end{remark}

\begin{remark}
We note that, in contrast to the other obtained results, it is not clear whether the result of \autoref{thm: est H' sigma and q} transfers to an analogous result for solutions to SDEs. This is because the proof of \autoref{thm: est H' sigma and q} relies on self-similarity and stationarity of increments of the process $X=\sigma Z^{H,q}$ and typical solutions to SDEs do not have these two properties. 
\end{remark}

\section{Auxiliary lemmas}
\label{app:aux_lemmas}

In this appendix, we collect some auxiliary lemmas that are used throughout the article. 

\begin{lemma}
\label{lem: pVar bounded by Lp}
Let $H\in (1/2,1)$ and $q\in\N$. Let also $f\in L^2(0,T)$ be a given deterministic function and let process $X=(X_t)_{t\in [0,T]}$ be defined by $X_t := \int_0^t f_s \d Z_s^{H,q}$, $t\in [0,T]$. Then
    \[\E V_{n,[0,T]}^2(X) \lesssim \delta_n^{2H - 1} \left\| f  \right\|_{L^2(0,T)}^2\]
\end{lemma}

\begin{proof}
By \autoref{lem:bddness_of_int} and the H\"older inequality, we have
    \begin{equation}
    \label{eq:L2norm_est_f}
    	\left\|  \int_0^T f_s \bm{1}_{[u,v]}(s)  \d Z_s^{H,q} \right\|_{L^2(\Omega)} \lesssim \left\| f \, \bm{1}_{[u,v]} \right\|_{L^{\frac{1}{H}}[0,T]} \leq (v-u)^{H-\frac{1}{2}}\left\| f \, \bm{1}_{[u,v]} \right\|_{L^2(0,T)}
    \end{equation}
for any $0 \leq u < v \leq T$. By appealing to this estimate, we can write
    \begin{multline*}
    \E V_{n,[0,T]}^2 (Y)  
    	= \sum_{i=0}^{n-1} \E \left|\int_0^T f_s \bm{1}_{[t_i^n,t_{i+1}^n]}(s) \d Z_s^{H,q} \right|^2 
    	\\
    	\lesssim  \sum_{i=0}^{n-1}  \delta_n^{2H-1}\left\| f \, \bm{1}_{[t_i^n,t_{i+1}^n]} \right\|_{L^2(0,T)}^2 
    	= \delta_n^{2H - 1} \left\| f  \right\|_{L^2(0,T)}^2. 
    \end{multline*}
\end{proof}

\begin{lemma}
\label{lem: difference of pVars}
    Let $f,g \in L^2(0,T)$ be two deterministic functions and consider processes $X=(X_t)_{t\in [0,T]}$, defined by $X_t := \int_0^t f_s \d Z_s^{H,q}$, and $Y=(Y_t)_{t\in [0,T]}$, defined by $Y_t := \int_0^t g_s \d Z_s^{H,q}$. Then
    \[\E |V_{n,[0,T]}^2(X) - V_{n,[0,T]}^2(Y)| \lesssim \delta_n^{2H - 1} \left\| f-g \right\|_{L^2(0,T)} \left(\left\| f \right\|_{L^2(0,T)} + \left\| g \right\|_{L^2(0,T)} \right)\]
\end{lemma}

\begin{proof}
We can proceed similarly as in the proof of \cite[Lemma 4.2]{GueNua05} with $1/H$ replaced by $2$ and get 
		\begin{equation*}
		\E |V_{n,[0,T]}^2(X) - V_{n,[0,T]}^2(Y)| \leq   2 (\E V_{n,[0,T]}^2(X-Y))^\frac{1}{2} ( (\E V_{n,[0,T]}^2(X))^{\frac{1}{2}} +(\E V_{n,[0,T]}^2(Y))^{\frac{1}{2}}).
	  \end{equation*}
Application of \autoref{lem: pVar bounded by Lp} concludes the proof.
\end{proof}

\begin{lemma}[Delta method for convergence in probability]
\label{lem:delta_method_for_p}
Let $\{X_n\}_n$ and $X$ be random variables defined on a common probability space, let $\{a_n\}_n\subseteq\R$ be such that $a_n \to \infty$, and let $\theta\in\R$. Let also $g$ be a real-valued function that is continuously differentiable in the neighbourhood of $\theta$. Assume that 
	\[ a_n (X_n-\theta) \xrightarrow[n\to\infty]{\P} X.\]
Then 
	\[ a_n (g(X_n)-g(\theta))\xrightarrow[n\to\infty]{\P} g'(\theta) X.\]
\end{lemma}

\begin{proof}
By using the mean value theorem, we have for $n\in\N$ that 
	\[ g(X_n) -g(\theta) = g'(\theta_n) (X_n-\theta), \quad a.s.\]
where $\theta_n$ is a random variable that lies between $X_n$ and $\theta$ almost surely. By the continuous mapping theorem, we then have $g'(\theta_n)\xrightarrow[n\to\infty]{\P} g'(\theta)$. It then follows that 
	\[ a_n(g(X_n)-g(\theta)) = g'(\theta_n) a_n(X_n-\theta) \xrightarrow[n\to\infty]{\P} g'(\theta)X\]
by the assumption of the lemma. 
\end{proof}

\end{document}